\documentclass{article}
\pdfoutput=1

\usepackage[narrow]{arxiv}

\usepackage[utf8]{inputenc} 
\usepackage[T1]{fontenc}    
\usepackage[hypertexnames=false]{hyperref}       
\usepackage{url}            
\usepackage{booktabs}       
\usepackage{amsfonts}       
\usepackage{nicefrac}       
\usepackage{microtype}      

\usepackage{lipsum}         
\usepackage{graphicx}
\usepackage{natbib}
\usepackage{doi}

\title{Non-Euclidean SGD\\ for Structured Optimization:\\ Unified Analysis and Improved Rates}
\renewcommand{\shorttitle}{Non-Euclidean SGD for Structured Optimization}

\newif\ifuniqueAffiliation
\uniqueAffiliationtrue

\ifuniqueAffiliation 
\author{
  Dmitry Kovalev \\
  Yandex Research \\
  \texttt{dakovalev1@gmail.com}
  \And
  Ekaterina Borodich\\
  Yandex Research \\
  \texttt{borodich.ed@phystech.edu}
}
\else
\usepackage{authblk}

\author[1]{%
  Dmitry Kovalev\thanks{\texttt{dakovalev1@gmail.com}}%
}
\affil[1]{Yandex Research}
\fi

\usepackage{algorithm,algpseudocodex}
\usepackage{nicematrix}
\usepackage{makecell}

\usepackage[textsize=tiny,color=orange!20,bordercolor=orange,]{todonotes}
\presetkeys{todonotes}{inline}{}

\usepackage{xcolor}
\hypersetup{
  colorlinks=true,
  linkcolor=red!75!black,
  citecolor=blue!50!black
}

\usepackage{mathtools,amsthm,amssymb}
\usepackage{bm}
\allowdisplaybreaks[4]

\usepackage{enumitem}
\setlist{topsep=0pt,partopsep=0pt,itemsep=0pt}

\usepackage{cleveref}
\usepackage{crossreftools}
\crefname{problem}{problem}{problems}
\crefname{condition}{condition}{conditions}
\crefname{property}{property}{properties}
\creflabelformat{problem}{(#2#1#3)}

\usepackage{annotate}
\usepackage{mythm}
\usepackage{mymath}

\newcommand{\Spp}{\S_{\succ}}
\newcommand{\Rp}{\R_{\geq}}

\newcommand{\normt}{\gennorm{\mathrm{tr}}{}}
\newcommand{\sqnt}{\gennorm{\mathrm{tr}}{2}}
\newcommand{\normo}{\gennorm{\mathrm{op}}{}}

\newcommand{\bg}{\mathrm{D}}
\newcommand{\N}{\mathbb{N}}
\newcommand{\ox}{\overline{x}}
\newcommand{\dist}{\cR\brr}
\newcommand{\sqdist}{\cR^2\brr}
\newcommand{\distd}{\cR_*\brr}
\newcommand{\sqdistd}{\cR_*^2\brr}
\renewcommand{\L}{\mathbb{L}}

\makeatletter
\newcommand{\grad}[1]{\cG\mymath@subscript{#1}\brr}
\makeatother

\begin{document}
\maketitle

\begin{abstract}
  Recently, several instances of non-Euclidean SGD, including SignSGD, Lion, and Muon, have attracted significant interest from the optimization community due to their practical success in training deep neural networks. Consequently, a number of works have attempted to explain this success by developing theoretical convergence analyses. Unfortunately, these results cannot properly justify the superior performance of these methods, as they could not beat the convergence rate of vanilla Euclidean SGD. We resolve this important open problem by developing a new unified convergence analysis under the structured smoothness and gradient noise assumption. In particular, our results indicate that non-Euclidean SGD {\bf (i)} can exploit the sparsity or low-rank structure of the upper bounds on the Hessian and gradient noise, {\bf (ii)} can provably benefit from popular algorithmic tools such as extrapolation or momentum variance reduction, and {\bf (iii)} can match the state-of-the-art convergence rates of adaptive and more complex optimization algorithms such as AdaGrad and Shampoo.

\end{abstract}

\section{Introduction}

Optimization algorithms play a fundamental role in deep neural network training. Over the past couple of decades, a large number of optimization methods have been developed and successfully applied in practice. These range from vanilla stochastic gradient descent (SGD) with momentum \citep{robbins1951stochastic,polyak1964some} to more advanced and complex adaptive algorithms such as AdaGrad \citep{duchi2011adaptive}, Shampoo \citep{gupta2018shampoo}, and Adam/AdamW \citep{kingma2014adam,loshchilov2017decoupled}. The latter have become the standard choice in most deep learning tasks.

However, recently, the algorithmic memory requirements have been tightening up due to the significant growth of model sizes, especially of the current large language models (LLMs) \citep{dubey2024llama,liu2024deepseek,achiam2023gpt,team2023gemini}. Consequently, the optimization community has started shifting its attention to more memory-efficient, non-Euclidean variants of SGD. These include SignSGD \citep{bernstein2018signsgd}, Lion \citep{chen2023symbolic}, and Muon \citep{jordan2024muon}, the latter starting to replace AdamW as the standard choice in training LLMs \citep{team2025kimi,zeng2025glm,liu2025muon}. Unfortunately, despite the practical effectiveness of these methods, their theoretical understanding remains highly limited. In particular, their known theoretical convergence rates cannot even surpass the rates of vanilla SGD, as we will discuss later. In this paper, we aim to resolve this important open question and advance the theoretical understanding of non-Euclidean SGD methods.




\subsection{Non-Euclidean SGD: Overview}

\textbf{Euclidean methods.}
The convergence properties of Euclidean methods are relatively well-studied. In particular, \citet{ghadimi2013stochastic} established the iteration complexity of vanilla SGD for non-convex Lipschitz-smooth functions with uniformly bounded gradient noise. Additionally, \citet{cutkosky2020momentum} analyzed the convergence of normalized SGD with momentum under the same assumptions. Both complexity results have optimal dependence on the target precision $\epsilon > 0$, i.e., they cannot be improved by any stochastic first-order method without extra assumptions \citep{arjevani2023lower}.


\textbf{Sign-based methods.} The most notable sign-based SGD methods include SignSGD with momentum \citep{bernstein2018signsgd} and Lion \citep{chen2023symbolic}. The convergence rates for these methods were established by \citet{sun2023momentum} and \citet{dong2024convergence}, respectively. Note that we will not consider Lion further because its theoretical properties are very similar to those of SignSGD with momentum. In particular, generalizing the results of this paper to Lion is straightforward.

\textbf{Orthogonolized methods.} Orthogonalized SGD methods are specifically designed for functions defined on the space of $m\times n$ matrices. They apply the so-called orthogonalization procedure $G \mapsto (GG^\top)^{\frac{\dagger}{2}}G$ to the current gradient estimator, where $G \in \R^{m\times n}$. Orthogonalization was first proposed by \citet{carlson2015preconditioned} and was later rediscovered by \citet{jordan2024muon}, leading to the development of Muon -- an orthogonalized SGD with momentum. This algorithm has remarkable practical performance across a wide range of deep learning tasks. \citet{pethick2025training,kovalev2025understanding} developed a convergence analysis of Muon under non-Euclidean Lipschitz smoothness with respect to the spectral norm and uniformly bounded variance with respect to the Euclidean norm. Furthermore, because of Muon's exceptional performance, many works have appeared studying various aspects of this method \citep{khaled2025muonbp,zhang2025adagrad,gruntkowska2025drop,crawshaw2025exploration,wang2025muon,shen2025convergence,he2025low,zhang2025provable,chen2025muon,he2025demuon,lau2025polargrad,huang2025limuon,riabinin2025gluon,shulgin2025beyond,liu2025mars,sato2025analysis}.

\textbf{General non-Euclidean methods.} Motivated by the interest in sign-based and orthogonalized gradient methods, \citet{pethick2025training,kovalev2025understanding} developed general analysis of non-Euclidean SGD with momentum and arbitrary norms, which recovers the convergence results for Muon and SignSGD.

\textbf{Connection with adaptive methods.}
Recently, \citet{xie2025structured,kovalev2025sgd} developed a unified convergence analysis for the adaptive meta-algorithm of \citet{gupta2017unified}, which recovers many popular adaptive methods as special cases, including AdaGrad-Norm \citep{streeter2010less}, AdaGrad, and ASGO/One-sided Shampoo \citep{an2025asgo,xie2025structured}, a more streamlined variant of Shampoo. More importantly, \citet{xie2025structured,kovalev2025sgd} showed how various instances of the adaptive meta-algorithm can exploit non-Euclidean smoothness and bounded gradient noise properties or their variants with respect to a family of non-Euclidean norms.

\subsection{Summary of Contributions}

\textbf{Main question.} As previously mentioned, the main question considered in this paper is that, to the best of our knowledge, {\em a proper theoretical explanation for the superior practical performance of non-Euclidean SGD is missing in the literature}.
For instance, existing theoretical results are unable to explain the practical success of Muon, since its state-of-the-art theoretical convergence guarantees are worse than those of vanilla SGD. More concretely, the state-of-the-art convergence rate of Muon under non-Euclidean Lipschitz smoothness with respect to the spectral norm and uniformly bounded variance with respect to the Euclidean norm \citep[Corollary~2]{kovalev2025sgd} is $\min\brf{m,n}$ times worse than the rate of vanilla SGD under the same assumptions \citep[Equation~2.17]{ghadimi2013stochastic}\footnote{Note that non-Euclidean Lipschitz smoothness with respect to the spectral norm implies Euclidean smoothness with the same Lipschitz constant.}. This factor may be huge, especially in current deep neural networks. Similar issues apply to SignSGD and non-Euclidean SGD with arbitrary norms, where the current state-of-the-art convergence rates \citep{pethick2025training,kovalev2025understanding} cannot surpass Euclidean SGD.

\textbf{Our solution.} In this paper, we provide a comprehensive answer to the above question and develop a unified convergence analysis of non-Euclidean SGD with momentum and weight decay (\Cref{alg}) for a family of non-Euclidean norms. \Cref{alg} recovers normalized SGD with momentum \citep{cutkosky2020momentum}, SignSGD with momentum \citep{bernstein2018signsgd}, and Muon \citep{jordan2024muon}. Our analysis implies {\em new state-of-the-art convergence guarantees for Muon and SignSGD, which beat Euclidean SGD and may finally explain the superior practical performance of these methods}. In addition, we establish the convergence guarantees for \Cref{alg} with extrapolation (option~\ref{eq:m2}) and momentum variance reduction (option~\ref{eq:m3}) tricks, which provably accelerate convergence. Finally, under the additional convexity assumption for the objective function, we show that {\em Muon and SignSGD with momentum variance reduction can match the existing state-of-the-art convergence rates for their respective adaptive counterparts: ASGO/One-sided Shampoo \citep{an2025asgo,xie2025structured} and AdaGrad \citep{duchi2011adaptive,liu2024adagrad}}. Thus, our results support the recent trend of abandoning adaptive methods in favor of simpler and more memory-efficient instances of non-Euclidean SGD.


\section{Preliminaries}

Formally, in this paper, we are concerned with the problem of minimizing a continuously differentiable objective function $f\brr{} \colon \cX \to \R$ defined on a finite-dimensional Euclidean space $\cX$.

\textbf{Notation.}
$\N_0$ is the set of non-negative integers;
$\Rp$ is the set of non-negative real numbers;
$\L$ is the space of linear operators $\cX \to \cX$;
$\S$ and $\Spp$ are the sets of self-adjoint and self-adjoint positive definite operators $\cX \to \cX$, respectively;
$\norm{x}_{\mB} = \sqrt{\<x,\mB x>}$, where $x \in \cX$ and $\mB \in \Spp$;
$\normt{\mB}$ and $\normo{\mB}$ are the nuclear and spectral norms of an operator $\mB \in \L$;
for all $z \in \cX$, by $z\<z,> \in \S$ we denote the rank-1 self-adjoint operator $x\mapsto z\<z,x>$;
$\mI \in \S$ denotes the identity operator. For a vector $x \in \R^d$, by $\brs{x}_{[i]}$ we denote its $i$-th component.

\subsection{Assumptions on the Objective Function}\label{sec:ass}

In this section, we describe the assumptions that we impose on the objective function $f\brr{}$. Throughout this work, we will always use \Cref{ass:L,ass:grad} for the analysis of \Cref{alg}, while \Cref{ass:T,ass:M} will only be used for the analysis of momentum with extrapolation and momentum variance reduction options, respectively. Note that all \Cref{ass:L,ass:grad,ass:T,ass:M} are given in terms of weighted Euclidean norms of the form $\norm{}_\mH$, where $\mH \in \Spp$. The main reason is that when combined with \Cref{ass:H}, they allow us to describe the non-Euclidean structure in the properties of the objective function $f\brr{}$. More details are given in \Cref{sec:structure}.

\textbf{Smoothness and unbiased gradient estimator.}
We assume that the objective function $f\brr{}$ is $1$-smooth with respect to the weighted Euclidean norm $\norm{}_\mL$, where $\mL \in \Spp$ is a self-adjoint positive definite operator. That is, the following \Cref{ass:L} holds:
\begin{aequation}<ass:L>
  \abs{\bg_f(x;x')} \leq \tfrac{1}{2}\sqn{x-x'}_\mL
  \quad\text{for all}\;\; x,x' \in \cX.
\end{aequation}
We also assume that there exists a stochastic estimator $g(\cdot;\xi) \colon \cX \to \cX$ of the gradient $\nabla f\brr{}$, where $\xi$ is a random variable sampled from some probability distribution $\cD$. We assume that $g(\cdot;\xi)$ is unbiased and has bounded variance, i.e., the following \Cref{ass:grad} holds:
\begin{aequation}<ass:grad>
  \E[\xi \sim \cD]{g(x;\xi)} = \nabla f(x)
  \quad\text{and}\quad
  \E[\xi \sim \cD]{\sqn{n(x;\xi)}_{\mSigma^{-1}}} \leq \normt{\mSigma}
  \quad\text{for all}\;\; x \in \cX,
\end{aequation}
where $\mSigma \in \Spp$ is a self-adjoint positive definite operator, $n(x;\xi) = g(x;\xi) - \nabla f(x)$ is the noise of the gradient estimator.

\textbf{Second-order smoothness for momentum with extrapolation.} To achieve improved convergence rates for SGD using momentum with extrapolation (\cref{eq:m2} in \Cref{alg}), we assume that the objective function $f(x)$ is second-order Lipschitz-smooth with respect to the weighted Euclidean norm $\norm{}_\mT$, where $\mT\in\Spp$. That is, the following \Cref{ass:T} holds:
\begin{aequation}<ass:T>
  \sqn{\brs{\nabla^2 f(x) - \nabla^2 f(x')}(x-x')}_{\mT^{-1}} \leq \normt{\mT}^{-1} \cdot \norm{x-x'}_\mT^4
  \quad\text{for all}\;\; x,x' \in \cX.
\end{aequation}
Note that the simplified version of this assumption with $\mT \propto \mI$ was previously used by \citet{cutkosky2020momentum} to obtain improved convergence rates for normalized SGD with momentum and extrapolation.

\textbf{Mean-square Lipschitz gradients for momentum variance reduction.}
To achieve even better convergence rates for SGD using momentum with extrapolation (\cref{eq:m2} in \Cref{alg}), we assume that the stochastic gradient estimator $g(\cdot;x)$ is mean-square $1$-Lipschitz with respect to the weighted Euclidean norm $\norm{}_\mM$, where $\mM \in \Spp$. That is, the following \Cref{ass:M} holds:
\begin{aequation}<ass:M>
  \E[\xi \sim \cD]{\sqn{g(x;\xi) - g(x';\xi)}_{\mM^{-1}}}  \leq \sqn{x-x'}_\mM
  \quad\text{for all}\;\; x,x' \in \cX.
\end{aequation}
Note that the simplified version of this assumption with $\mM \propto \mI$ and its variants was previously used by \citet{cutkosky2019momentum,fatkhullin2022sharp} to  analyze SGD with momentum variance reduction.

\textbf{On the restrictiveness of the assumptions.}
In the case where $\mL, \mSigma \propto \mI$, \Cref{ass:L,ass:grad} reduces to the standard Lipschitz continuity of the gradient $\nabla f\brr{}$ and the uniformly bounded noise of the gradient estimator $g(\cdot;\xi)$, which may not be the best choice for modeling the properties of complex functions such as deep neural networks. For instance, \citet{zhang2019gradient,chen2023generalized,gorbunov2024methods,vankov2024optimizing,borodich2025nesterov} analyzed the convergence of gradient methods under the $(L_0,L_1)$-smoothness, and \citet{vaswani2019fast,faw2022power,attia2023sgd,sadchikov2024local} analyzed SGD under the affine variance. Both are more realistic assumptions than Lipschitz smoothness and bounded variance. However, we do not consider these settings because our main goal is to show how non-Euclidean SGD can exploit the structure in the operators $\mL, \mSigma$ (and $\mT,\mM$ in \Cref{ass:T,ass:M}), and applying the techniques from our paper to generalized variants of \Cref{ass:L,ass:grad} should be straightforward.

\subsection{Main Structural Assumption}\label{sec:structure}

\begin{table}[t]
  \caption{Instances of \Cref{alg} induced by the choice of the linear space $\cX$, the space of operators $\cH$ satisfying \Cref{ass:H}, and the non-Euclidean norm $\dist{}$ defined in \cref{eq:dist}: normalized SGD \citep{cutkosky2020momentum}, SignSGD \citep{bernstein2018signsgd}, and Muon \citep{jordan2024muon}, along with their adaptive counterparts: AdaGrad-Norm \citep{streeter2010less}, AdaGrad \citep{duchi2011adaptive}, and ASGO/One-sided Shampoo \citep{an2025asgo,xie2025structured}.}
  \begin{center}
    \begin{NiceTabular}[]{ccccc}
      \toprule
      \bf \Cref{alg} & \bf Adaptive Counterpart &$\cX$& $\cH$ & $\dist{}$
      \\\midrule
      Normalized SGD & AdaGrad-Norm &$\R^d$& $\brf{g \mapsto \beta g : \beta \in \R}$ & $\frac{1}{\sqrt{d}}\norm{}_2$
      \\\midrule
      SignSGD & AdaGrad &$\R^d$& $\brf{g  \mapsto \bm{b} \odot g : \bm{b} \in \R^d} $ & $\normi{}$
      \\\midrule
      Muon & One-sided Shampoo &$\R^{m \times n}$& $\brf{G \mapsto B G : B \in \S^m}$ & $\frac{1}{\sqrt{n}}\normo{}$
      \\\bottomrule
    \end{NiceTabular}
  \end{center}
  \label{tab:H}
\end{table}

We impose the additional structural assumptions on the operators $\mL,\mM,\mT,\mSigma \in \Spp$. In particular, we assume that $\mL$, $\mM$, $\mT$, and $\mSigma$ belong to a certain space of linear operators $\cH \subset \S$ that satisfies the conditions in \Cref{ass:H}.
\begin{assumption}\label{ass:H}
  $\mL,\mM,\mT,\mSigma \in \cH\cap\Spp$, where $\cH \subset \S$ is a linear subspace of self-adjoint operators, satisfying the following conditions:
  \begin{enumerate}[label=\bf(\roman{*})]
    \item The space $\cH$ contains identity operator: $\mI \in \cH$.
    \item The space $\cH$ is closed under symmetric multiplications: $\mH\mF\mH \in \cH$ for all $\mH,\mF \in \cH$
    \item The orthogonal projection onto $\cH$ is order-preserving: $\proj[\cH]{\mH} \in \Spp$ for all $\mH \in \Spp$.
  \end{enumerate}
\end{assumption}
Next, we define functions $\dist{}\colon \cX \to \Rp$ and $\distd{}\colon \cX \to \Rp$ as follows:
\begin{equation}\label{eq:dist}
  \dist{x} = \normo{\sqrt{\proj[\cH]{\mX}}}
  \quad\text{and}\quad
  \distd{x} =  \normt{\sqrt{\proj[\cH]{\mX}}},
  \quad\text{where}\quad
  \mX = x\<x,>.
\end{equation}
One can verify that both functions $\dist{}$ and $\distd{}$ are norms, possibly non-Euclidean, as indicated by \Cref{lem:dist}. Moreover, \Cref{lem:dist} implies that these norms are dual to each other.
\begin{lemma}<lem:dist>
  $\dist{}$ and $\distd{}$ are norms, that is, they are subadditive, absolutely homogeneous, and positive-definite functions. Moreover, these norms are dual to each other:
  \begin{equation}
    \dist{x} = \tsup_{\distd{y} \leq 1}\<x,y>
    \quad\text{and}\quad
    \distd{y} = \tsup_{\dist{x} \leq 1}\<x,y>.
  \end{equation}
\end{lemma}

In \Cref{tab:H}, we describe the main instances of \Cref{alg} induced by the choice of the operator space $\cH$ and the corresponding non-Euclidean norm $\dist{}$: normalized SGD and SignSGD with momentum and weight decay, and Muon with weight decay. In addition, we provide details on how these algorithms exploit the structure of the operators $\mL,\mM,\mT,\mSigma$ in \Cref{ass:H} for improved convergence guarantees in \Cref{sec:discussion}.

\textbf{Connection with non-Euclidean setting.}
In the following \Cref{lem:HR}, we establish the connection between the non-Euclidean norm $\dist{}$ and an arbitrary weighted Euclidean norm of the form $\norm{}_\mH$, where $\mH \in \cH \cap \Spp$.
\begin{lemma}<lem:HR>
  Let \Cref{ass:H} hold and let $\mH \in \cH \cap \Spp$. Then the following inequalities hold:
  \begin{equation}
    \sqn{x}_\mH \leq \normt{\mH}\cdot\sqdist{x}
    \quad\text{and}\quad
    \sqn{x}_{\mH^{-1}} \geq \normt{\mH}^{-1}\cdot\sqdistd{x}
    \quad\text{for all}\;\; x \in \cX.
  \end{equation}
\end{lemma}

Using \Cref{lem:HR}, it is not hard to verify that the assumptions on the objective function from \Cref{sec:ass} imply their non-Euclidean counterparts, that is, for all $x,x' \in \cX$, the following relations hold:
\begin{align}
  \label{nE:L}
  \text{\Cref{ass:L}}
  \quad&\Rightarrow\quad
  &
  \abs{\bg_f(x;x')} &\leq \tfrac{1}{2}\normt{\mL}\cdot\dist{x-x'}
  ;\\
  \label{nE:Sigma}
  \text{\Cref{ass:grad}}
  \quad&\Rightarrow\quad
  &
  \E*[\xi\sim\cD]{\sqdistd{n(x;\xi)}} &\leq \sqnt{\mSigma}
  ;\\
  \label{nE:T}
  \text{\Cref{ass:T}}
  \quad&\Rightarrow\quad
  &
  \distd*{\brs{\nabla^2 f(x) - \nabla^2 f(x')}(x-x')} &\leq \normt{\mT}\cdot\sqdist{x-x'}
  ;\\
  \label{nE:M}
  \text{\Cref{ass:M}}
  \quad&\Rightarrow\quad
  &
  \E*[\xi\sim\cD]{\sqdistd{g(x;\xi) - g(x';\xi)}} &\leq \sqnt{\mM}\cdot\sqdist{x-x'}.
\end{align}
It is worth highlighting that the reverse implications do not hold. Indeed, \Cref{ass:L,ass:grad,ass:T,ass:M} combined with \Cref{ass:H} do not imply the non-Euclidean counterparts in general. However, the additional structure implied by these assumptions is exactly what allows us to achieve the main goal: justifying the efficiency of non-Euclidean SGD. This also contrasts with the existing state-of-the-art results \citep{pethick2025training,kovalev2025understanding}, which used the non-Euclidean assumptions above, but could not improve upon the baseline Euclidean SGD.

\textbf{Connection with adaptive methods.}
\citet{xie2025structured,kovalev2025sgd} developed a unified convergence analysis for the adaptive-meta algorithm \citep{gupta2017unified}, which recovers AdaGrad-Norm, AdaGrad, and ASGO/One-sided Shampoo as special cases. As we will see further, these algorithms are closely related to normalized SGD, SignSGD, and Muon, respectively. In particular, our \Cref{ass:H} is similar to Assumption~1 by \citet{kovalev2025sgd}, and the definition of $\dist{}$ in \cref{eq:dist} is identical to that by \citet{kovalev2025sgd}, used in the analysis of the adaptive meta-algorithm. More details are available in \Cref{sec:tr}.

\subsection{Constrained Optimization Setting}\label{sec:Q}

Further, in this paper, we focus on the following constrained variant of the optimization problem:
\begin{equation}\label[problem]{eq:mainQ}
  f^* = \min_{x \in Q} f(x),
\end{equation}
where $Q = \brf{x \in \cX: \dist{x} \leq \cR}$ is a ball of possibly infinite radius $\cR \in [0,+\infty]$, defined in terms of the non-Euclidean norm $\dist{}$. As we will see further, this additional constraint does not introduce any overhead, as it is automatically handled by adding weight decay in \Cref{alg}. On the other hand, we can solve the unconstrained version of the problem by simply choosing a large enough radius $\cR$. Note that the presence of this constraint is often desirable: it can be seen as a regularization, which may encourage dense \citep{chen2023lion} or high-rank \citep{chen2025muon} solution in the case of sign-based or orthogonalized methods, respectively.

\textbf{Convergence criterion.}
We define the function $\grad{\cR}{}\colon Q \to \Rp$, which we are going to use as a convergence criterion, as follows:
\begin{equation}\label{eq:criterion}
  \grad{\cR}{x} =
  \begin{cases}
    \tfrac{1}{\cR}\max_{ \dist{x'} \leq \cR} \<\nabla f(x), x-x'> & \cR \in [0,+\infty)\\
    \distd{\nabla f(x)} & \cR = +\infty
  \end{cases}.
\end{equation}
In the case $\cR = +\infty$, the function $\grad{\cR}{}$ is simply the dual gradient norm $\distd{\nabla f\brr{}}$, which is the standard convergence criterion for the unconstrained non-convex setting. In the constrained case, $\cR < +\infty$, the function $\grad{\cR}{}$ coincides with the Frank-Wolfe gap \citep{jaggi2013revisiting} up to the scaling factor $1/\cR$. Note that the scaling factor $1/\cR$ ensures that the function $\cR \mapsto \grad{\cR}{x}$ is continuous at $\cR = +\infty$ for all $x \in \cX$. In addition, $\grad{\cR}{}$ coincides with the convergence criterion used by \citet{dong2024convergence,chen2025muon}. For completeness, in \Cref{lem:criterion}, we verify that the function $\grad{\cR}{}$ is indeed a meaningful convergence criterion in the constrained case: null points of $\grad{\cR}{}$ are the first-order stationary points for \cref{eq:mainQ}. Additionally, in the constrained convex case, $\grad{\cR}{}$ can be used to upper-bound the objective function optimality gap, as stated in \Cref{lem:convex}.

\begin{lemma}<lem:criterion>
  $\grad{\cR}{x^*} = 0$ if and only if $x^* \in Q$ is a first-order stationary point in \cref{eq:mainQ}, that is, $-\nabla f(x^*) \in \cN_Q(x^*)$, where $\cN_Q(x) = \brf{y \in \cX : \<y,x-x'> \geq 0 \text{\;for all\;}x' \in Q}$ is the normal cone to the set $Q$ at the point $x \in Q$.
\end{lemma}

\begin{lemma}<lem:convex>
  Let $\cR < +\infty$ and let the objective function $f(x)$ be convex. Then the following inequality holds for all $x \in Q$:
  \begin{equation}
    \cR \cdot \grad{\cR}{x} \geq f(x) - f^*.
  \end{equation}
\end{lemma}

\section{Non-Euclidean SGD with Momentum and Weight Decay}\label{sec:theory}

\begin{algorithm}[t]
  \caption{Non-Euclidean SGD with Momentum and Weight Decay}
  \label{alg}
  \begin{algorithmic}[1]
    \State \textbf{input:} $x_0 \in Q$
    \label{line:input}
    \State \textbf{parameters:} momentum $\alpha\in (0,1)$, weight decay $\beta \in [0,1)$, stepsize $\eta > 0$
    \State sample $\xi_0\sim \cD$ and compute $m_0 = g(x_0;\xi_0)$
    \label{line:init}
    \For{$k=0,1,\dots,K-1$}
    \State compute $x_{k+1}$ as follows:
    \begin{align}\label{eq:x}
      x_{k+1} = \argmin_{x \in \cX}\;\<m_k,x-x_k>
      \quad\text{s.t.}\quad \dist{x-(1-\beta)x_k}\leq \eta
    \end{align}
    \State sample $\xi_{k+1}\sim \cD$ and compute $m_{k+1}$ using one of the following options:
    \label{line:m}
    \begin{subequations}
      \begin{align}
        \label{eq:m1}
        m_{k+1} &= (1-\alpha)m_k + \alpha g(x_{k+1},\xi_{k+1})
        \\
        \label{eq:m2}
        m_{k+1} &= (1-\alpha)m_k + \alpha g(\ox_{k+1},\xi_{k+1})\\
        \label{eq:m3}
        m_{k+1} &= (1-\alpha)(m_k - g(x_k,\xi_{k+1})) + g(x_{k+1};\xi_{k+1})
      \end{align}
    \end{subequations}
    where $\ox_{k+1} = x_k + (1/\alpha)(x_{k+1} - x_k)$.
    \EndFor
    \State \textbf{output:} $x_K \in \cX$
  \end{algorithmic}
\end{algorithm}

\subsection{Non-Euclidean Trust-Region Gradient Step}\label{sec:tr}

In this section, we develop the convergence analysis for the non-Euclidean SGD with momentum and weight decay, which is formalized as \Cref{alg}.
We start with the non-Euclidean trust-region gradient step with the gradient estimator $m_{k+1}$ in \cref{eq:x}. It was analyzed in details for general non-Euclidean norms \Cref{eq:dist} by \citet{kovalev2025understanding}. Further, we develop the analysis under the structural assumptions in \Cref{sec:structure}. We start with the basic properties in \Cref{lem:x}.
\begin{lemma}<lem:x>
  Let \Cref{ass:H} hold and let $\beta = \eta/\cR$. Then, $\beta\dist{x_k} \leq \eta$ and $\dist{x_{k+1} - x_k} \leq 2\eta$.
\end{lemma}
Next, we define the operator $\mB \in \Spp$ as follows:
\begin{equation}\label{eq:B}
  \mB = \tfrac{1}{4}\normt{\mL}^{-1}\mL + \tfrac{1}{4}\normt{\mM}^{-1}\mM + \tfrac{1}{4}\normt{\mSigma}^{-1}\mSigma + \tfrac{1}{4}\normt{\mT}^{-1}\mT.
\end{equation}
It is easy to verify that $\mB \in \cH$ due to \Cref{ass:H}, and that $\normt{\mB} = 1$. We also provide useful bounds on the operator $\mB$ in \Cref{lem:B} and the main result for the trust-region step in \Cref{lem:tr}.
\begin{lemma}<lem:B>
  Let \Cref{ass:H} hold. Then for $\mH \in \brf{\mL,\mM,\mSigma,\mT}$, the following inequalities hold:
  \begin{equation}
    \mH \preceq 4\normt{\mH}\cdot \mB,\quad
    4\normt{\mH}\cdot\mH^{-1} \succeq \mB^{-1}
  \end{equation}
\end{lemma}

\begin{lemma}<lem:tr>
  Let \Cref{ass:L,ass:H} hold and let $\beta= \eta/\cR$. Then the following inequality holds, where we use the convention $1/0 = +\infty$:
  \begin{equation}
    \eta \grad{\eta/\beta}{x_k} \leq f(x_k) - f(x_{k+1}) + 4\eta\norm{\nabla f(x_k) - m_k}_{\mB^{-1}} + 2\normt{\mL}\eta^2.
  \end{equation}
\end{lemma}

\textbf{Connection with adaptive methods.}
As previously mentioned, \citet{xie2025structured,kovalev2025sgd} developed a unified analysis of the adaptive meta-algorithm. More concretely, they analyzed the preconditioned SGD of the form $x_{k+1} = x_k - \mH_k g_k$, where $g_k \in \cX$ is the current gradient estimator. Under \Cref{ass:H}, the preconditioning operator can be expressed as $\mH_k = \eta \brs{\delta\mI+ \proj[\cH]{\mS_k}}^{-1/2}$, where $\mS_k = \tsum_{i=0}^{k}g_i\<g_i,>$ is the cumulative sum of gradient outer squares, and $\eta,\delta > 0$ are positive parameters \citep{kovalev2025sgd}. One can verify that when the gradient accumulation is turned off, i.e., $\mS_k = g_k\<g_k,>$, the preconditioned SGD step above is equivalent to the trust-region step in \cref{eq:x} with $m_k = g_k$ and $\beta = 0$. Hence, the adaptive meta-algorithm can be seen as an adaptive counterpart to \Cref{alg}. Consequently, AdaGrad and ASGO/One-sided Shampoo are adaptive ``twins'' to SignSGD and Muon, as indicated in \Cref{tab:H}.


\begin{table}[t]
  \caption{Summary of the results in \Cref{sec:theory}. The number of iterations $K$ of \Cref{alg} required to ensure $\E{\min_{k\in\brf{0,\dots,K-1}}\grad{\cR}{x_k}} \leq \epsilon$ in the non-convex case and $\E{f(x_K) - f^*}\leq \epsilon$ in the convex case is given. Only the leading terms (with the worst dependence on $\epsilon$) are reported. Universal and logarithmic factors are omitted. Extra assumptions beyond \Cref{ass:L,ass:grad,ass:H} are provided.}

  \begin{center}
    \begin{NiceTabular}[cell-space-limits=0.3em]{cllc}
      \toprule
      \bf Reference
      &
      \bf Non-Convex Case
      &
      \bf Convex Case
      &
      \bf Extra Assumption
      \\\hline
      \Cref{thm:m1}
      &
      $
      \displaystyle
      \tmfrac{\sqnt{\mSigma}\normt{\mL}\Delta_0}{\epsilon^4}
      $
      &
      $
      \displaystyle
      \tmfrac{\sqnt{\mSigma}\normt{\mL}\cR^4}{\epsilon^3}
      $
      &
      ---
      \\\hline
      \Cref{thm:m2}
      &
      $
      \displaystyle
      \tmfrac{\sqnt{\mSigma}\normt{\mT}^{0.5}\Delta_0}{\epsilon^{3.5}}
      $
      &
      $
      \displaystyle
      \tmfrac{\sqnt{\mSigma}\normt{\mT}^{0.5}\cR^{3.5}}{\epsilon^{2.5}}
      $
      &
      \Cref{ass:T}
      \\\hline
      \Cref{thm:m3}
      &
      $
      \displaystyle
      \tmfrac{\normt{\mSigma}^3 + \normt{\mSigma}\normt{\mM}\Delta_0}{\epsilon^3}
      $
      &
      $
      \displaystyle
      \tmfrac{\sqnt{\mSigma}\cR^2 + \normt{\mSigma}\normt{\mM}\cR^3}{\epsilon^2}
      $
      &
      \Cref{ass:M}
      \\\bottomrule
    \end{NiceTabular}
  \end{center}
  \label{tab:rates}
\end{table}

\subsection{Main Convergence Results}\label{sec:main}

In this section, we obtain the main convergence results for \Cref{alg}, which are summarized in \Cref{tab:rates}. First, in the following \Cref{lem:m1,,lem:m2,,lem:m3}, we analyze three options for computing the gradient estimator $m_{k+1}$ on \cref{line:m} using momentum, momentum with extrapolation, and momentum variance reduction in \cref{eq:m1,,eq:m2,,eq:m3}, respectively.
\begin{lemma}<lem:m1>
  Let \Cref{ass:grad,,ass:H} hold and let $\beta= \eta/\cR$. Let $m_{k+1}$ be defined using \cref{eq:m1} at each iteration of \Cref{alg}. Then the following inequality holds:
  \begin{equation}
    \E{\sqn{m_{k+1} - \nabla f(x_{k+1})}_{\mB^{-1}}}
    \leq
    (1-\alpha)\E{\sqn{m_k - \nabla f(x_k)}_{\mB^{-1}}}
    +64\sqnt{\mL}\eta^2/\alpha
    +4\alpha^2\sqnt{\mSigma}.
  \end{equation}
\end{lemma}
\begin{lemma}<lem:m2>
  Let \Cref{ass:T,,ass:grad,,ass:H} hold and let $\beta= \eta/\cR$. Let $m_{k+1}$ be defined using \cref{eq:m2} at each iteration of \Cref{alg}. Then the following inequality holds:
  \begin{equation}
    \E{\sqn{m_{k+1} - \nabla f(x_{k+1})}_{\mB^{-1}}}
    \leq
    (1-\alpha)\E{\sqn{m_k - \nabla f(x_k)}_{\mB^{-1}}}
    +192\sqnt{\mT}\eta^4/\alpha^3
    +4\alpha^2\sqnt{\mSigma}.
  \end{equation}
\end{lemma}
\begin{lemma}<lem:m3>
  Let \Cref{ass:M,,ass:grad,,ass:H} hold and let $\beta = \eta/\cR$. Let $m_{k+1}$ be defined using \cref{eq:m3} at each iteration of \Cref{alg}. Then the following inequality holds:
  \begin{equation}
    \E{\sqn{m_{k+1} - \nabla f(x_{k+1})}_{\mB^{-1}}}
    \leq
    (1-\alpha)^2\E{\sqn{m_k - \nabla f(x_k)}_{\mB^{-1}}}
    +128\sqnt{\mM}\eta^2
    +8\alpha^2\sqnt{\mSigma}.
  \end{equation}
\end{lemma}

Using the above lemmas, we  obtain  the following \Cref{cor:m1,,cor:m2,,cor:m3}

\begin{corollary}<cor:m1>
  Under the conditions of \Cref{lem:m1}, the following inequality holds:
  \begin{equation}
    \E{\norm{m_{k} - \nabla f(x_{k})}_{\mB^{-1}}}
    \leq
    2\brs*{(1-\alpha/2)^k + \sqrt{\alpha}}\normt{\mSigma}
    +8(\eta/\alpha)\normt{\mL}.
  \end{equation}
\end{corollary}
\begin{corollary}<cor:m2>
  Under the conditions of \Cref{lem:m2}, the following inequality holds:
  \begin{equation}
    \E{\norm{m_{k} - \nabla f(x_{k})}_{\mB^{-1}}}
    \leq
    2\brs*{(1-\alpha/2)^k + \sqrt{\alpha}}\normt{\mSigma}
    +8\sqrt{3}(\eta/\alpha)^2\normt{\mT}.
  \end{equation}
\end{corollary}
\begin{corollary}<cor:m3>
  Under the conditions of \Cref{lem:m3}, the following inequality holds:
  \begin{equation}
    \E{\norm{m_{k} - \nabla f(x_{k})}_{\mB^{-1}}}
    \leq
    2\brs{(1-\alpha)^k + \sqrt{2\alpha}}\normt{\mSigma}
    +8\sqrt{2/\alpha}\normt{\mM}\eta.
  \end{equation}
\end{corollary}
Now, we are ready to establish the convergence rates for \Cref{alg} with momentum, momentum with extrapolation, and momentum variance reduction in the following \Cref{thm:m1,,thm:m2,,thm:m3}, respectively. In particular, we get the number of iterations $K$ that are required to reach the precision $\E*{\min_{k\in\brf{0,K-1}} \grad{\cR}{x_k}} \leq \epsilon$ in the non-convex case, and the precision $\E{f(x_K) - f^*} \leq \epsilon$ in the convex case. Note that we also cover the non-convex unconstrained case, where $\cR = +\infty$. Therefore, \Cref{thm:m1,,thm:m2,,thm:m3} treat the constrained and unconstrained settings in a unified manner.

\begin{theorem}<thm:m1>
  Let \Cref{ass:L,,ass:grad,,ass:H} hold.
  Let $m_{k+1}$ be defined using \cref{eq:m1} at each iteration of \Cref{alg}.
  Then, choosing the parameters of \Cref{alg} as
  \begin{equation}\label{eq:parameters1}
    \alpha = \min\brf*{1,\cO\brr*{\tmfrac{\epsilon^2}{\sqnt{\mSigma}}}}
    ,\quad
    \beta = \min\brf*{1,\cO\brr*{\tmfrac{\epsilon\alpha}{\normt{\mL}\cR}}}
    ,\quad
    \eta = \beta\cR,
  \end{equation}
  and the number of iterations $K \in \N$ of \Cref{alg} as
  \begin{equation}\label{eq:K1}
    K = \cO\brr*{\max\brf*{
        1,
        \tmfrac{\Delta_0}{\epsilon \cR},
        \tmfrac{\normt{\mL}\Delta_0}{\epsilon^2},
        \tmfrac{\normt{\mSigma}^3}{\epsilon^3},
        \tmfrac{\sqnt{\mSigma}\normt{\mL}\Delta_0}{\epsilon^4}
    }},
  \end{equation}
  where $\Delta_0 = f(x_0) - f^*$, implies $\E*{\min_{k\in\brf{0,K-1}} \grad{\cR}{x_k}} \leq \epsilon$.
  Furthermore, if the function $f(x)$ is assumed to be convex, then choosing the parameters of \Cref{alg} as
  \begin{equation}\label{eq:parameters_cvx1}
    \alpha = \min\brf*{1,\cO\brr*{\tmfrac{\epsilon^2}{\sqnt{\mSigma}\cR^2}}}
    ,\quad
    \beta = \min\brf*{\alpha,
      \cO\brr*{
        \tmfrac{\epsilon\alpha}{\normt{\mL}\cR^2}
    }}
    ,\quad
    \eta = \beta\cR,
  \end{equation}
  and the number of iterations $K \in \N$ of \Cref{alg} as
  \begin{equation}\label{eq:K_cvx1}
    K = \tilde\cO\brr*{\max\brf*{
        1,
        \tmfrac{\normt{\mL}\cR^2}{\epsilon},
        \tmfrac{\sqnt{\mSigma}\cR^2}{\epsilon^2},
        \tmfrac{\sqnt{\mSigma}\normt{\mL}\cR^4}{\epsilon^3}
    }}.
  \end{equation}
  implies $\E{f(x_K) - f^*} \leq \epsilon$.
\end{theorem}

\begin{theorem}<thm:m2>
  Let \Cref{ass:L,,ass:grad,,ass:T,,ass:H} hold.
  Let $m_{k+1}$ be defined using \cref{eq:m2} at each iteration of \Cref{alg}.
  Then, choosing the parameters of \Cref{alg} as
  \begin{equation}\label{eq:parameters2}
    \alpha = \min\brf*{1,\cO\brr*{\tmfrac{\epsilon^2}{\sqnt{\mSigma}}}}
    ,\quad
    \beta = \min\brf*{1,\cO\brr*{\tmfrac{\epsilon}{\normt{\mL}\cR}},\cO\brr*{\tmfrac{\sqrt{\epsilon}\alpha}{\sqrt{\normt{\mT}}\cR}}}, \quad
    \eta = \beta\cR.
  \end{equation}
  and the number of iterations $K \in \N$ of \Cref{alg} as
  \begin{equation}\label{eq:K2}
    K = \cO\brr*{\max\brf*{
        1,
        \tmfrac{\Delta_0}{\epsilon\cR},
        \tmfrac{\normt{\mL}\Delta_0}{\epsilon^2},
        \tmfrac{\normt{\mSigma}^3}{\epsilon^3},
        \tmfrac{\sqnt{\mSigma}\normt{\mT}^{0.5}\Delta_0}{\epsilon^{3.5}}
    }},
  \end{equation}
  where $\Delta_0 = f(x_0) - f^*$, implies $\E*{\min_{k\in\brf{0,K-1}} \grad{\cR}{x_k}} \leq \epsilon$.
  Furthermore, if the function $f(x)$ is assumed to be convex, then choosing the parameters of \Cref{alg} as
  \begin{equation}\label{eq:parameters_cvx2}
    \alpha = \min\brf*{1,\cO\brr*{\tmfrac{\epsilon^2}{\sqnt{\mSigma}}}}
    ,\quad
    \beta = \min\brf*{\alpha,\cO\brr*{\tmfrac{\epsilon}{\normt{\mL}\cR}},\cO\brr*{\tmfrac{\sqrt{\epsilon}\alpha}{\textstyle\sqrt{\normt{\mT}\cR^3}}}}, \quad
    \eta = \beta\cR.
  \end{equation}
  and the number of iterations $K \in \N$ of \Cref{alg} as
  \begin{equation}\label{eq:K_cvx2}
    K = \tilde\cO\brr*{\max\brf*{
        1,
        \tmfrac{\normt{\mL}\cR^2}{\epsilon},
        \tmfrac{\sqnt{\mSigma}\cR^2}{\epsilon^2},
        \tmfrac{\sqnt{\mSigma}\normt{\mT}^{0.5}\cR^{3.5}}{\epsilon^{2.5}}
    }}.
  \end{equation}
  implies $\E{f(x_K) - f^*} \leq \epsilon$.
\end{theorem}

\begin{theorem}<thm:m3>
  Let \Cref{ass:L,,ass:grad,,ass:M,,ass:H} hold.
  Let $m_{k+1}$ be defined using \cref{eq:m3} at each iteration of \Cref{alg}.
  Then, choosing the parameters of \Cref{alg} as
  \begin{equation}\label{eq:parameters3}
    \alpha = \min\brf*{1,\cO\brr*{\tmfrac{\epsilon^2}{\sqnt{\mSigma}}}}
    ,\quad
    \beta = \min\brf*{1,\cO\brr*{\tmfrac{\epsilon}{\normt{\mL}\cR}},\cO\brr*{\tmfrac{\epsilon\sqrt{\alpha}}{\normt{\mM}\cR}}}, \quad
    \eta = \beta\cR,
  \end{equation}
  and the number of iterations $K \in \N$ of \Cref{alg} as
  \begin{equation}\label{eq:K3}
    K = \cO\brr*{\max\brf*{
        1,
        \tmfrac{\Delta_0}{\epsilon\cR},
        \tmfrac{\normt{\mL}\Delta_0}{\epsilon^2},
        \tmfrac{\normt{\mM}\Delta_0}{\epsilon^2},
        \tmfrac{\normt{\mSigma}^3}{\epsilon^3},
        \tmfrac{\normt{\mSigma}\normt{\mM}\Delta_0}{\epsilon^3}
    }},
  \end{equation}
  where $\Delta_0 = f(x_0) - f^*$, implies $\E*{\min_{k\in\brf{0,K-1}} \grad{\cR}{x_k}} \leq \epsilon$.
  Furthermore, if the function $f(x)$ is assumed to be convex, then choosing the parameters of \Cref{alg} as
  \begin{equation}\label{eq:parameters_cvx3}
    \alpha = \min\brf*{1,\cO\brr*{\tmfrac{\epsilon^2}{\sqnt{\mSigma}\cR^2}}}
    ,\quad
    \beta = \min\brf*{
      \alpha,
      \cO\brr*{\tmfrac{\epsilon}{\normt{\mL}\cR^2}},
      \cO\brr*{\tmfrac{\epsilon\sqrt{\alpha}}{\normt{\mM}\cR^2}}
    }
    ,\quad
    \eta = \beta\cR,
  \end{equation}
  and the number of iterations $K \in \N$ of \Cref{alg} as
  \begin{equation}\label{eq:K_cvx3}
    K = \tilde\cO\brr*{\max\brf*{
        1,
        \tmfrac{\normt{\mL}\cR^2}{\epsilon},
        \tmfrac{\sqnt{\mSigma}\cR^2}{\epsilon^2},
        \tmfrac{\normt{\mSigma}\normt{\mM}\cR^3}{\epsilon^2}
    }}.
  \end{equation}
  implies $\E{f(x_K) - f^*} \leq \epsilon$.
\end{theorem}

\subsection{Discussion}\label{sec:discussion}

\begin{table}[t]
  \caption{Comparison of the iterations complexities in the convex case for non-Euclidean SGD with momentum variance reduction (\Cref{alg}, \Cref{thm:m3}) with the results for adaptive meta-algorithm \citep{kovalev2025sgd,xie2025structured}. Universal constant factors are omitted.}

  \begin{center}
    \begin{NiceTabular}[cell-space-limits=0.3em]{lcc}
      \toprule
      \bf Method & \Cref{alg} (Option \ref{eq:m3}) & Adaptive Meta-Algorithm
      \\\hline
      \bf Complexity
      &
      $
      \displaystyle
      \tmfrac{\normt{\mL}\cR^2}{\epsilon}
      +
      \tmfrac{\sqnt{\mSigma}\cR^2}{\epsilon^2}
      +
      \tmfrac{\normt{\mSigma}\normt{\mM}\cR^3}{\epsilon^2}
      $
      &
      $
      \displaystyle
      \tmfrac{\normt{\mL}\cR^2}{\epsilon}
      +
      \tmfrac{\sqnt{\mSigma}\cR^2}{\epsilon^2}
      $
      \\\hline
      \bf Assumptions & \labelcref{ass:L,ass:grad,ass:M,ass:H} & \labelcref{ass:L,ass:grad,ass:H}
      \\
      \bottomrule
    \end{NiceTabular}
  \end{center}
  \label{tab:adaptive}
\end{table}

First, we compare the convergence results for Non-Euclidean SGD with momentum variance reduction (\Cref{alg}, Option~\ref{eq:m3}) in \Cref{thm:m3} with the convergence results for the adaptive meta-algorithm by \citet[Theorems~1 and~3]{kovalev2025sgd} in the convex setting. This comparison is summarized in \Cref{tab:adaptive}. \Cref{thm:m3} requires the additional \Cref{ass:M}, which is slightly more restrictive. Additionally, both rates have optimal dependence on $\epsilon$ \citep{nemirovskij1983problem,lan2012optimal} and match as long as $\normt{\mM}\cR \lesssim \normt{\mSigma}$.

Next, we show that \Cref{alg} can exploit the main structural \Cref{ass:H}. In particular, we show that sign-based and orthogonalized variants of \Cref{alg} can provably outperform the Euclidean variant of \Cref{alg}. For simplicity, we consider Option~\ref{eq:m1}, although the comparison below clearly applies to Options~\ref{eq:m2} and~\ref{eq:m3} as well.

\textbf{SignSGD vs Normalized SGD.} In the case of sign-based methods, the properties of the objective function are measured in terms of the weighted Euclidean norms of the form $\sqn{x}_{\bm{h}} = \<x, \bm{h} \odot x>$, where $\bm{h} \in \R_+^d$ and $x \in \cX = \R^d$. Consequently, \Cref{ass:L,ass:grad} reduce to
\begin{equation}
  \abs{\bg_f(x;x')} \leq \tfrac{1}{2}\sqn{x-x'}_{ \diag{\bm{l}}}
  ,\qquad
  \E[\xi \sim \cD]{ \brs{n_\xi (x)}_{[i]}^2} \leq \brs{\bm{\sigma}}_{[i]}^2 \text{\;for all\;}i,
\end{equation}
where $\bm{l}, \bm{\sigma} \in \R_+^d$. Note that, for simplicity and for illustrative purposes, we use here a slightly more restrictive variant of \Cref{ass:grad}, namely, the coordinate-wise bound on the variance. Under these assumptions, the convergence rates in \Cref{tab:rates} (\Cref{thm:m1}, convex case) reduce to
\begin{equation}
  \text{SignSGD:}\quad
  \tmfrac{\sqn{\bm{\sigma}}_1\norm{\bm{l}}_1 \normi{x^*}^4}{\epsilon^3}
  ,\qquad
  \text{Normalized SGD:}\quad
  \tmfrac{\sqn{\bm{\sigma}}_2\normi{\bm{l}} \norm{x^*}_2^4}{\epsilon^3}.
\end{equation}
We can conclude that the complexity of SignSGD improves upon the complexity of Normalized SGD, when the solution $x^* \in \R^d$ is dense, i.e., $\norm{x^*}_2 \sim \sqrt{d} \cdot \normi{x^*}$, and when the coordinates in the objective function have very different scales, i.e., the vectors $\bm{l},\bm{\sigma}$ have a few dominant components that are much greater than the rest, which implies $\norm{\bm{l}}_1 \ll d \cdot \normi{\bm{l}}$, $\norm{\bm{\sigma}}_1 \ll \sqrt{d} \cdot \norm{\bm{\sigma}}_2$. This observation aligns with and complements the existing theory for sign-based methods \citep{liu2024adagrad,jiang2024convergence}.

\textbf{Muon vs Normalized SGD.}
In the case of orthogonalized methods, the properties of the objective function are measured in terms of the weighted Euclidean norms of the form $\sqn{X}_{H} = \<X, H X>$, where $H \in \Spp^{m}$ and $X \in \cX = \R^{m\times n}$. \Cref{ass:L,ass:grad} reduce to the following:
\begin{equation}
  \abs{\bg_f(X;X')} \leq \tfrac{1}{2}\sqn{X-X'}_{L}
  ,\qquad
  \E[\xi\sim\cD]{N_\xi(X) N_\xi(X)^\top} \preceq \Sigma^2,
\end{equation}
where $L,\Sigma \in \Spp^m$, and $N_\xi(X) \in \R^{m\times n}$ is the gradient noise reshaped into an $m \times n$ matrix. Under these assumptions, the convergence rates in \Cref{tab:rates} (\Cref{thm:m1}, convex case) reduce to
\begin{equation}
  \text{Muon:}\quad
  \tmfrac{\sqnt{\Sigma}\normt{L} \normo{X^*}^4}{\epsilon^3}
  ,\qquad
  \text{Normalized SGD:}\quad
  \tmfrac{\sqn{\Sigma}_F\normo{L} \norm{X^*}_F^4}{\epsilon^3}.
\end{equation}
We can conclude that the complexity of Muon improves upon the complexity of Normalized SGD, when the solution $X^* \in \R^{m\times n}$ is high-rank, i.e., $\norm{X^*}_F \sim \sqrt{r}\cdot \normo{X^*}$, where $r = \min\brf{m,n}$, and when the matrices $L,\Sigma$ are approximately low-rank, i.e., $\normt{L} \ll r\cdot \normo{L}$, $\normt{\Sigma} \ll \sqrt{r}\cdot \norm{\Sigma}_F$. This observation aligns with and complements the existing theory for orthogonalized/matrix-based optimizers \citep{an2025asgo}.

\newpage

\cite*{}

\bibliographystyle{apalike}
\bibliography{references}

\newpage
\appendix
\part*{Appendix}

\section[Proofs from \crtcref{sec:structure}]{Proofs from \Cref{sec:structure}}

\proofsubsection{lem:dist}

\textbf{Part I.}
From \Cref{ass:H} it follows that $\mH^k \in \cH$ for all $\mH \in \cH$, $k \in \N_0$. Based on this fact, using polynomial interpolation, it is easy to show that $\cH$ is closed under arbitrary operator functions, that is, $\psi(\mH) \in \cH$ for all $\mH \in \cH$ and $\psi\brr{}\colon \R\to\R$. Hence, we can use Lemma~4 of \citet{kovalev2025sgd} to prove that the function $\dist{}$ is a norm.

\textbf{Part II.}
We start by proving the following Fenchel-Young inequality for all $x,y\in\cX$:
\begin{equation}\label{eq:FY}
  \<x,y> \leq \tfrac{1}{2}\sqdist{x} + \tfrac{1}{2}\sqdist{y}.
\end{equation}
Indeed, we can lower-bound $\tfrac{1}{2}\sqdist{x}$ as follows:
\begin{align*}
  \tfrac{1}{2}\sqdist{x}
  &\at{uses the definition of $\dist{}$ in \cref{eq:dist} and the properties of $\normo{}$}=
  \<x,y> + \sup_{\normt{\mH}\leq 1,\mH \in \Spp} \brs*{\tfrac{1}{2}\<\proj[\cH]{\mX},\mH> - \<x,y>}
  \\&\geq
  \<x,y> + \adjustlimits\inf_{\;x\in \cX\;}\sup_{\normt{\mH}\leq 1,\mH \in \Spp} \brs*{\tfrac{1}{2}\<\proj[\cH]{\mX},\mH> - \<x,y>}
  \\&\at{uses the fact that $\inf_{x}\sup_\mH \brs{\dots} \geq \sup_\mH\inf_{x} \brs{\dots}$}\geq
  \<x,y> + \adjustlimits\sup_{\normt{\mH}\leq 1,\mH \in \Spp}\inf_{x\in \cX} \brs*{\tfrac{1}{2}\<\proj[\cH]{\mX},\mH> - \<x,y>}
  \\&\at{uses the properties of the projection and the definition of $\mX$ in \cref{eq:dist}}\geq
  \<x,y> + \adjustlimits\sup_{\normt{\mH}\leq 1,\mH \in \Spp}\inf_{x\in \cX} \brs*{\tfrac{1}{2}\sqn{x}_{\proj[\cH]{\mH}} - \<x,y>}
  \\&\at{uses the fact that $\proj[\cH]{\mH} \in \Spp$ due to \Cref{ass:H}}=
  \<x,y> - \inf_{\normt{\mH}\leq 1,\mH \in \Spp} \tfrac{1}{2}\sqn{y}_{\proj[\cH]{\mH}^{-1}},
\end{align*}
where \annotate.
After rearranging, we obtain the following inequality for all $\mH \in \Spp$ such that $\normt{\mH} \leq 1$:
\begin{equation}
  \<x,y> \leq \tfrac{1}{2}\sqdist{x} + \tfrac{1}{2}\sqn{y}_{\proj[\cH]{\mH}^{-1}}.
\end{equation}
Let $\delta > 0$ and let $\mH \in \S$ be defined as follows:
\begin{equation}
  \mH = \normt{\sqrt{\delta \mI + \proj[\cH]{\mY}}}^{-1}\sqrt{\delta \mI + \proj[\cH]{\mY}},
  \quad\text{where}\quad \mY = y\<y,>.
\end{equation}
Using \Cref{ass:H} and the fact that $\cH$ is closed under operator functions from Part~I, we can show that $\mH,\mH^{-1} \in \Spp \cap \cH$. Hence, we can upper-bound $\<x,y>$ as follows:
\begin{align*}
  \<x,y>
  &\at{uses the upper bound on $\<x,y>$ above}\leq
  \tfrac{1}{2}\sqdist{x} + \tfrac{1}{2}\sqn{y}_{\proj[\cH]{\mH}^{-1}}
  \\&\at{uses the fact that $\mH \in \cH$}=
  \tfrac{1}{2}\sqdist{x} + \tfrac{1}{2}\sqn{y}_{\mH^{-1}}
  \\&\at{uses the definition of $\mY$ above and the fact that $\mH \in \Spp$}\leq
  \tfrac{1}{2}\sqdist{x} + \tfrac{1}{2}\<\delta \mI + \mY,\mH^{-1}>
  \\&\at{uses the properties of the projection and the fact that $\mH^{-1} \in \cH$}=
  \tfrac{1}{2}\sqdist{x} + \tfrac{1}{2}\<\delta \mI + \proj[\cH]{\mY} ,\mH^{-1}>
  \\&\at{uses the definition of $\mH$ above}=
  \tfrac{1}{2}\sqdist{x} + \tfrac{1}{2}\sqnt{\sqrt{\delta \mI + \proj[\cH]{\mY}}},
\end{align*}
where \annotate. After taking the limit $\delta \to +0$ and using the definition of $\distd{}$ in \cref{eq:dist}, we obtain the desired inequality in \cref{eq:FY}.

Next, it is easy to show that $\distd{}$ is absolutely homogeneous. Hence, using \cref{eq:FY}, is easy to show that $\<x,y> \leq \distd{y}$ for all $x \in \cX$ such that $\dist{x} \leq 1$. This immediately implies the inequality $\tsup_{\dist{x} \leq 1}\<x,y> \leq \distd{y}$. Now, we will prove the reverse inequality. Let $\delta > 0$ and let $x \in \cX$ be defined as follows:
\begin{equation}
  x = \mF y,
  \quad\text{where}\quad
  \mF = \brs{\delta\mI + \proj[\cH]{\mY}}^{-1/2}.
\end{equation}
Using \Cref{ass:H} and the fact that $\cH$ is closed under operator functions from Part~I, we can show that $\mF \in \Spp \cap \cH$. Hence, we can verify that $\dist{x} \leq 1$ as follows:
\begin{align*}
  \sqdist{x}
  &\at{uses the definition of $\dist{}$ in \cref{eq:dist} and the properties of $\normo{}$}=
  \sup_{\normt{\mH}\leq 1,\mH \in \Spp} \<\proj[\cH]{\mX},\mH>
  \\&\at{uses the properties of the projection}=
  \sup_{\normt{\mH}\leq 1,\mH \in \Spp} \<\mX,\proj[\cH]{\mH}>
  \\&\at{uses the definitions of $x$ and $\mY$ above and the definition of $\mX$ in \cref{eq:dist}}=
  \sup_{\normt{\mH}\leq 1,\mH \in \Spp} \<\mF\mY\mF,\proj[\cH]{\mH}>
  \\&=
  \sup_{\normt{\mH}\leq 1,\mH \in \Spp} \<\mY,\mF\proj[\cH]{\mH}\mF>
  \\&\at{uses the properties of the projection and the fact that $\mF\proj[\cH]{\mH}\mF \in \cH$ due to \Cref{ass:H}}=
  \sup_{\normt{\mH}\leq 1,\mH \in \Spp} \<\proj[\cH]{\mY},\mF\proj[\cH]{\mH}\mF>
  \\&\at{uses the fact that $\proj[\cH]{\mH} \in \Spp$ due to \Cref{ass:H}}\leq
  \sup_{\normt{\mH}\leq 1,\mH \in \Spp} \<\delta \mI + \proj[\cH]{\mY},\mF\proj[\cH]{\mH}\mF>
  \\&\at{uses the definition of $\mF$ above}=
  \sup_{\normt{\mH}\leq 1,\mH \in \Spp} \<\mF^{-2},\mF\proj[\cH]{\mH}\mF>
  \\&\at{uses the properties of the projection and fact that $\mI \in \cH$ due to \Cref{ass:H}}=
  \sup_{\normt{\mH}\leq 1,\mH \in \Spp} \<\mI,\mH> = 1,
\end{align*}
where \annotate. Hence, we can lower-bound $\tsup_{\dist{x} \leq 1}\<x,y>$ as follows:
\begin{align*}
  \tsup_{\dist{x} \leq 1}\<x,y>
  &\at{uses the definitions of $x$ and $\mY$ above, the properties of the projection, and the fact that $\mF \in \cH$}\geq
  \tsup_{\delta > 0}\<\mF,\proj[\cH]{\mY}>
  \\&\at{uses the definition of $\mF$ above}=
  \tsup_{\delta > 0}\<\brs{\delta\mI + \proj[\cH]{\mY}}^{-1/2},\proj[\cH]{\mY}>
  \\&=
  \tsup_{\delta > 0} \brs{\normt{\brs{\delta\mI + \proj[\cH]{\mY}}^{1/2}}- \delta\normt{\brs{\delta\mI + \proj[\cH]{\mY}}^{-1/2}}}
  \\&\at{is obtained by taking the limit $\delta \to +0$}\geq
  \normt{\brs{\proj[\cH]{\mY}}^{1/2}}
  \at{uses the definition of $\distd{}$ in \cref{eq:dist}}=
  \distd{y},
\end{align*}
where \annotate. Therefore, we can conclude that $\tsup_{\dist{x} \leq 1}\<x,y> = \distd{y}$.

\textbf{Pert III.} From Part~II, it follows that $\distd{}$ is a dual norm to $\dist{}$. This implies that $\dist{}$ is a dual norm to $\distd{}$, that is, $\tsup_{\distd{y} \leq 1}\<x,y> = \dist{x}$.\qed

\proofsubsection{lem:HR}
We can upper-bound $\sqn{x}_\mH$ as follows:
\begin{align*}
  \sqn{x}_\mH
  &\at{uses the definition of $\mX$ in \cref{eq:dist}}=
  \<\mH,\mX>
  \at{uses the assumption $\mH \in \cH$ and the properties of the projection}=
  \<\mH,\proj[\cH]{\mX}>
  \at{uses H\"older's inequality for Schatten norms}\leq
  \normt{\mH} \cdot \normo{\proj[\cH]{\mX}}
  \at{uses the definition of $\distd{}$ in \cref{eq:dist}}=
  \normt{\mH} \cdot \sqdist{x},
\end{align*}
where \annotate. Furthermore, we can lower-bound $\sqdistd{x}$ as follows:
\begin{align*}
  \sqdistd{x}
  &\at{uses \Cref{lem:dist}}=
  \tsup_{\dist{x'} \leq 1} \<x,x'>^2
  \at{uses the previously obtained inequality $\sqn{x'}_\mH \leq \normt{\mH}\cdot\sqdist{x'}$}\geq
  \tsup_{\sqn{x'}_\mH \leq \normt{\mH}} \<x,x'>^2
  =
  \normt{\mH} \cdot \sqn{x}_{\mH^{-1}},
\end{align*}
where \annotate.\qed

\newpage

\section[Proofs for \crtCref{sec:Q}]{Proofs for \Cref{sec:Q}}

\proofsubsection{lem:criterion}

Assume $\cR < +\infty$, otherwise the proof is trivial. Suppose that $x^* \in Q$ is a first-order stationary point for the optimization problem $\min_{x \in Q} f(x)$, that is, $-\nabla f(x^*) \in \cN_Q(x^*)$ By the definition of the normal cone, we have
\begin{align*}
  \<-\nabla f(x^*),x^*-x> \geq 0
  \quad \text{for all}\;\; x \in Q,
\end{align*}
which implies $\grad{\cR}{x^*} \leq 0$. Hence, using the fact that $\grad{\cR}{x} \geq 0$ for all $x \in Q$, we get $\grad{\cR}{x^*} = 0$.

Conversely, suppose that $\grad{\cR}{x^*} = 0$. Then, by the definition of $\grad{\cR}{}$ in \cref{eq:criterion}, we have
\begin{align*}
  \<\nabla f(x^*),x^* - x> \leq 0
  \quad \text{for all}\;\; x \in Q,
\end{align*}
which by the definition of the normal cone implies $-\nabla f(x^*) \in \cN_Q(x^*)$.\qed

\proofsubsection{lem:convex}

We can upper-bound $f(x) - f^*$ for all $x \in Q$ as follows:
\begin{align*}
  f(x) - f^*
  &\at{uses the definition of $f^*$ in \cref{eq:mainQ}}\leq
  \max_{x' \in Q} \brs{f(x) - f(x')}
  \at{uses the convexity assumption}\leq
  \max_{x' \in Q} \<\nabla f(x),x - x'>
  \at{uses the definition of $\grad{\cR}{x}$ in \cref{eq:criterion} and the definition of $Q$}\leq
  \cR \cdot \grad{\cR}{x},
\end{align*}
where \annotate.\qed

\newpage
\section[Proofs for \crtCref{sec:tr}]{Proofs for \Cref{sec:tr}}

\proofsubsection{lem:x}

We prove the inequality $\beta \dist{x_k} \leq \eta$ by induction. The base case $k=0$ is implied by the inclusion $x_0 \in Q$ on \cref{line:input} and the definition of $Q$. Furthermore, we assume the inequality $\beta \dist{x_k} \leq \eta$ and prove the inequality $\beta \dist{x_{k+1}} \leq \eta$ as follows:
\begin{align*}
  \beta \dist{x_{k+1}}
  &\at{uses the properties of the norm $\dist{}$}\leq
  \beta \dist{x_{k+1} - (1-\beta)x_k} + (1-\beta) \beta\dist{x_k}
  \\&\at{uses the definition of $x_{k+1}$ in \cref{eq:x}}\leq
  \beta \eta + (1-\beta) \beta\dist{x_k}
  \at{uses the induction hypothesis}\leq
  \beta \eta + (1-\beta) \eta
  = \eta,
\end{align*}
where \annotate.
Next, we prove the inequality $\dist{x_{k+1} - x_k} \leq 2\eta$ as follows:
\begin{align*}
  \dist{x_{k+1} - x_k}
  &\at{uses the properties of the norm $\dist{}$}\leq
  \dist{x_{k+1} - (1-\beta)x_k}
  +\beta\dist{x_k}
  \\&\at{uses the definition of $x_{k+1}$ in \cref{eq:x}}\leq
  \eta
  +\beta\dist{x_k}
  \at{uses the previously obtained inequality above}\leq
  2\eta,
\end{align*}
where \annotate.\qed

\proofsubsection{lem:B}

The inequality $\mH \preceq 4\normt{\mH}\cdot \mB$ is trivially implied by the definition of $\mB$ in \cref{eq:B} and \Cref{ass:H}. Hence, using the L\"owner-Heinz theorem \citep[Theorem~2.6]{carlen2010trace}, we obtain the inequality $4\normt{\mH}\cdot\mH^{-1} \succeq \mB^{-1}$.\qed

\proofsubsection{lem:tr}

In the case $\beta > 0$, we can upper-bound $f(x_{k+1})$ as follows:
\begin{align*}
  f(x_{k+1})
  &\at{uses \Cref{ass:L,ass:H} and \Cref{lem:HR}}\leq
  f(x_k) + \<\nabla f(x_k), x_{k+1} - x_k > + \tfrac{1}{2}\normt{\mL}\sqdist{x_{k+1} - x_k}
  \\&\at{uses \Cref{lem:x}}\leq
  f(x_k) + \<\nabla f(x_k) - m_k, x_{k+1} - x_k > + \<m_k, x_{k+1} - x_k > + 2\normt{\mL}\eta^2
  \\&\at{uses the definition of $x_{k+1}$ in \cref{eq:x}}\leq
  f(x_k) + \<\nabla f(x_k) - m_k, x_{k+1} - x_k > + \min_{\dist{x - (1-\beta)x_k}\leq \eta}\<m_k, x - x_k > + 2\normt{\mL}\eta^2
  \\&\at{uses the variable exchange $x = \beta z + (1-\beta)x_k$}=
  f(x_k) + \<\nabla f(x_k) - m_k, x_{k+1} - x_k > + \min_{\beta\dist{z}\leq \eta}\beta\<m_k, z - x_k > + 2\normt{\mL}\eta^2
  \\&\at{uses the definition of $\grad{\eta/\beta}{}$ in \cref{eq:criterion} and $z = \argmin_{\beta\dist{z}\leq\eta}\<\nabla f(x_k),z-x_k>$}\leq
  f(x_k) + \<\nabla f(x_k) - m_k, x_{k+1} - \beta z - (1-\beta)x_k> + 2\normt{\mL}\eta^2
  -\eta \grad{\eta/\beta}{x_k}
  \\&\at{uses the Cauchy-Schwarz inequality}\leq
  f(x_k) + \norm{x_{k+1} - \beta z - (1-\beta)x_k}_{\mB}\norm{\nabla f(x_k) - m_k}_{\mB^{-1}} + 2\normt{\mL}\eta^2
  -\eta \grad{\eta/\beta}{x_k}
  \\&\at{uses \Cref{lem:HR} and the definition of $\mB$ in \cref{eq:B}}\leq
  f(x_k) + \dist{x_{k+1} - \beta z - (1-\beta)x_k}\norm{\nabla f(x_k) - m_k}_{\mB^{-1}} + 2\normt{\mL}\eta^2
  -\eta \grad{\eta/\beta}{x_k}
  \\&\at{uses the triangle inequality, \Cref{lem:x} and the fact that $\beta\dist{z}\leq \eta$}\leq
  f(x_k) + 4\eta\norm{\nabla f(x_k) - m_k}_{\mB^{-1}} + 2\normt{\mL}\eta^2
  -\eta \grad{\eta/\beta}{x_k},
\end{align*}
where \annotate.
In the case $\beta = 0$, we can upper-bound $f(x_{k+1})$ as follows:
\begin{align*}
  f(x_{k+1})
  &\at{uses \Cref{ass:L,ass:H} and \Cref{lem:HR}}\leq
  f(x_k) + \<\nabla f(x_k), x_{k+1} - x_k > + \tfrac{1}{2}\normt{\mL}\sqdist{x_{k+1} - x_k}
  \\&\at{uses \Cref{lem:x}}\leq
  f(x_k) + \<\nabla f(x_k) - m_k, x_{k+1} - x_k > + \<m_k, x_{k+1} - x_k > + 2\normt{\mL}\eta^2
  \\&\at{uses the definition of $x_{k+1}$ in \cref{eq:x} and $\beta = 0$}\leq
  f(x_k) + \<\nabla f(x_k) - m_k, x_{k+1} - x_k > + \min_{\dist{x - x_k}\leq \eta}\<m_k, x - x_k > + 2\normt{\mL}\eta^2
  \\&\at{uses \Cref{lem:dist}}=
  f(x_k) + \<\nabla f(x_k) - m_k, x_{k+1} - x_k >  - \eta \distd{m_k} + 2\normt{\mL}\eta^2
  \\&\at{uses the Cauchy-Schwarz inequality}\leq
  f(x_k) + \norm{x_{k+1} - x_k}_{\mB}\norm{f(x_k) - m_k}_{\mB^{-1}} - \eta \distd{m_k} + 2\normt{\mL}\eta^2
  \\&\at{uses the triangle inequality}\leq
  f(x_k) + \norm{x_{k+1} - x_k}_{\mB}\norm{f(x_k) - m_k}_{\mB^{-1}} + \eta \distd{\nabla f(x_k) - m_k}
  - \eta \distd{\nabla f(x_k)} \\&+ 2\normt{\mL}\eta^2
  \\&\at{uses \Cref{lem:HR} and the definition of $\mB$ in \cref{eq:B}}\leq
  f(x_k) + \brs{\dist{x_{k+1} - x_k} + \eta}\norm{f(x_k) - m_k}_{\mB^{-1}}
  - \eta \distd{\nabla f(x_k)} + 2\normt{\mL}\eta^2
  \\&\at{uses \Cref{lem:x} and the inequality $3\leq 4$}\leq
  f(x_k) + 4\eta\norm{f(x_k) - m_k}_{\mB^{-1}}
  - \eta \distd{\nabla f(x_k)} + 2\normt{\mL}\eta^2
  \\&\at{uses the definition of $\grad{+\infty}{}$ in \cref{eq:criterion} and $\eta/\beta = +\infty$}=
  f(x_k) + 4\eta\norm{f(x_k) - m_k}_{\mB^{-1}} - \eta \grad{\eta/\beta}{x_k} + 2\normt{\mL}\eta^2,
\end{align*}
where \annotate.\qed

\newpage
\section[Proofs for \crtCref{sec:main} (Option~\ref{eq:m1})]{Proofs for \Cref{sec:main} (Option~\ref{eq:m1})}

\proofsubsection{lem:m1}

We can upper-bound $\E[\xi_{k+1}\sim \cD]{\sqn{m_{k+1} - \nabla f(x_{k+1})}_{\mB^{-1}}}$ as follows:
\begin{align*}
  \mi{1}\E[\xi_{k+1}\sim \cD]{\sqn{m_{k+1} - \nabla f(x_{k+1})}_{\mB^{-1}}}=
  \\&=
  \E[\xi_{k+1}\sim \cD]{\sqn{\E[\xi_{k+1}\sim \cD]{m_{k+1}} - \nabla f(x_{k+1})}_{\mB^{-1}}}
  +
  \E[\xi_{k+1}\sim \cD]{\sqn{m_{k+1} - \E[\xi_{k+1}\sim \cD]{m_{k+1}}}_{\mB^{-1}}}
  \\&\at{uses \cref{eq:m1} and the definition of $n(x;\xi)$ in \Cref{ass:grad}}\leq
  \sqn{(1-\alpha)(m_k- \nabla f(x_{k+1}))}_{\mB^{-1}}
  +
  \alpha^2\E[\xi_{k+1}\sim \cD]{\sqn{n(x_{k+1};\xi_{k+1})}_{\mB^{-1}}}
  \\&\at{uses \Cref{lem:B} and \Cref{ass:grad}}\leq
  \sqn{(1-\alpha)(m_k- \nabla f(x_k)) + (1-\alpha)(\nabla f(x_k)- \nabla f(x_{k+1}))}_{\mB^{-1}}
  +
  4\alpha^2\sqnt{\mSigma}
  \\&\at{uses the convexity of $\sqn{}_{\mB^{-1}}$}\leq
  (1-\alpha)\sqn{m_k - \nabla f(x_k)}_{\mB^{-1}}
  +((1-\alpha)^2/\alpha)\sqn{\nabla f(x_k) - \nabla f(x_{k+1})}_{\mB^{-1}}
  +4\alpha^2\sqnt{\mSigma}
  \\&\at{is implied by \Cref{ass:L} and \Cref{lem:B}}\leq
  (1-\alpha)\sqn{m_k - \nabla f(x_k)}_{\mB^{-1}}
  +(16\sqnt{\mL}/\alpha)\sqn{x_k - x_{k+1}}_{\mB}
  +4\alpha^2\sqnt{\mSigma}
  \\&\at{uses \Cref{lem:HR}, the definition of $\mB$ in \cref{eq:B}, and \Cref{lem:x}}\leq
  (1-\alpha)\sqn{m_k - \nabla f(x_k)}_{\mB^{-1}}
  +64\sqnt{\mL}\eta^2/\alpha
  +4\alpha^2\sqnt{\mSigma},
\end{align*}
where \annotate. It remains to take the full expectation.\qed

\proofsubsection{cor:m1}
We can upper-bound $\E{\sqn{m_{k+1} - \nabla f(x_{k+1})}_{\mB^{-1}}}$ as follows:
\begin{align*}
  \mi{3}\E*{\sqn{m_{k+1} - \nabla f(x_{k+1})}_{\mB^{-1}}}\leq
  \\&\at{uses \Cref{lem:m1}}\leq
  (1-\alpha)^{k+1}\E*{\sqn{m_0 - \nabla f(x_0)}_{\mB^{-1}}}
  +\brs*{64\sqnt{\mL}\eta^2/\alpha + 4\alpha^2\sqnt{\mSigma}} \tsum_{i = 0}^k (1 - \alpha)^i
  \\&\at{uses the definition $m_0 = g(x_0; \xi_0)$ on \cref{line:init} and the definition of $n(x;\xi)$ in \Cref{ass:grad}}\leq
  (1-\alpha)^{k+1}\E*{\sqn{n(x_0;\xi_0)}_{\mB^{-1}}}
  +\brs*{64\sqnt{\mL}\eta^2/\alpha + 4\alpha^2\sqnt{\mSigma}} \tsum_{i = 0}^k (1 - \alpha)^i
  \\&\at{uses \Cref{lem:B} and \Cref{ass:grad}}\leq
  4(1-\alpha)^{k+1}\sqnt{\mSigma}
  +\brs*{64\sqnt{\mL}\eta^2/\alpha + 4\alpha^2\sqnt{\mSigma}} \tsum_{i = 0}^k (1 - \alpha)^i
  \\&\at{uses the inequality $\tsum_{i = 0}^k (1 - \alpha)^i \leq 1/\alpha$}\leq
  4\brs*{(1-\alpha)^{k+1} + \alpha}\sqnt{\mSigma}
  +64(\eta/\alpha)^2\sqnt{\mL},
\end{align*}
where \annotate.
Now, we can upper-bound $\E*{\norm{m_k - \nabla f(x_k)}_{\mB^{-1}}}$ as follows:
\begin{align*}
  \E*{\norm{m_k - \nabla f(x_k)}_{\mB^{-1}}}
  &\at{uses Jensen's inequality for expectations}\leq
  \sqrt{\E*{\sqn{m_k - \nabla f(x_k)}_{\mB^{-1}}}}
  \\&\at{uses the upper bound above}\leq
  \sqrt{4\brs*{(1-\alpha)^k + \alpha}\sqnt{\mSigma}
  +64(\eta/\alpha)^2\sqnt{\mL}}
  \\&\at{uses the inequality $\sqrt{a+b}\leq\sqrt{a}+\sqrt{b}$ and the inequality $\sqrt{1-\alpha} \leq 1-\alpha/2$}\leq
  2\brs*{(1-\alpha/2)^k + \sqrt{\alpha}}\normt{\mSigma}
  +8\normt{\mL}\eta/\alpha,
\end{align*}
where \annotate.\qed

\proofsubsection{thm:m1}

We can upper-bound $\E*{\min_{k\in\brf{0,K-1}} \grad{\cR}{x_k}}$ as follows:
\begin{align*}
  \mi{3}
  \E*{\tmin_{k\in\brf{0,K-1}} \grad{\cR}{x_k}}\leq
  \\&\leq
  \tfrac{1}{K}\tsum_{k=0}^{K-1}\E*{ \grad{\cR}{x_k}}
  \\&\at{uses \Cref{lem:tr} and the definition of $\eta$ in \cref{eq:parameters1}}\leq
  \tfrac{1}{\eta K}\tsum_{k=0}^{K-1}\E*{ f(x_k) - f(x_{k+1}) + 4\eta\norm{\nabla f(x_k) - m_k}_{\mB^{-1}} + 2\normt{\mL}\eta^2}
  \\&\at{uses the definition of $\Delta_0$}\leq
  \tfrac{\Delta_0}{\eta K}
  +2\normt{\mL}\eta
  +\tfrac{4}{K}\tsum_{k=0}^{K-1}\E*{\norm{\nabla f(x_k) - m_k}_{\mB^{-1}}}
  \\&\at{uses \Cref{cor:m1}}\leq
  \tfrac{\Delta_0}{\eta K}
  +2\normt{\mL}\eta
  +\tfrac{4}{K}\tsum_{k=0}^{K-1}\brs*{2\brs*{(1-\alpha/2)^k + \sqrt{\alpha}}\normt{\mSigma}
  +8\normt{\mL}\eta/\alpha}
  \\&\at{uses the inequality $\tsum_{k=0}^{K-1}(1-\alpha/2)^k \leq 2/\alpha$}\leq
  \tfrac{\Delta_0}{\eta K}
  +\tfrac{16\normt{\mSigma}}{\alpha K}
  +\tfrac{32\normt{\mL}\eta}{\alpha}
  +2\normt{\mL}\eta
  +8\sqrt{\alpha}\normt{\mSigma},
\end{align*}
where \annotate. It remains to use the definitons of the parameters in \cref{eq:parameters1,eq:K1}.

Furthermore, in the convex case, we define the following Lyapunov function:
\begin{equation}\label{eq:Psi1}
  \Psi_k = f(x_k) - f^* + c\sqn{m_k - \nabla f(x_k)},
\end{equation}
where $c > 0$. We can upper-bound $\E{\Psi_{k+1}}$ as follows:
\begin{align*}
  \E{\Psi_{k+1}}
  &\at{uses \Cref{lem:tr,lem:m1}}\leq
  \E{f(x_k) - f^*} - \eta \E{\grad{\eta/\beta}{x_k}} + 4\eta\E{\norm{\nabla f(x_k) - m_k}_{\mB^{-1}}} + 2\normt{\mL}\eta^2
  \\&
  +c(1-\alpha)\E{\sqn{m_k - \nabla f(x_k)}_{\mB^{-1}}}
  +\tfrac{64c\eta^2}{\alpha}\sqnt{\mL}
  +4c\alpha^2\sqnt{\mSigma}
  \\&\at{uses the definition of $\eta$ in \cref{eq:parameters1} and \Cref{lem:convex}}\leq
  (1-\beta)\E{f(x_k) - f^*}
  +c(1-\alpha)\E{\sqn{m_k - \nabla f(x_k)}_{\mB^{-1}}}
  \\&
  +4\eta\E{\norm{\nabla f(x_k) - m_k}_{\mB^{-1}}}
  +2\normt{\mL}\eta^2
  +\tfrac{64c\eta^2}{\alpha}\sqnt{\mL}
  +4c\alpha^2\sqnt{\mSigma}
  \\&\at{uses Young's inequality}\leq
  (1-\beta)\E{f(x_k) - f^*}
  +c(1-\alpha/2)\E{\sqn{m_k - \nabla f(x_k)}_{\mB^{-1}}}
  +\tfrac{8\eta^2}{c\alpha}
  \\&
  +2\normt{\mL}\eta^2
  +\tfrac{64c\eta^2}{\alpha}\sqnt{\mL}
  +4c\alpha^2\sqnt{\mSigma}
  \\&\at{uses the definition of $\Psi_k$ in \cref{eq:Psi1}}\leq
  \brr*{1-\tfrac{1}{2}\min\brf{\alpha,\beta}}\E{\Psi_k}
  +2\normt{\mL}\eta^2
  +\tfrac{64c\eta^2}{\alpha}\sqnt{\mL}
  +4c\alpha^2\sqnt{\mSigma}
  +\tfrac{8\eta^2}{c\alpha}
  \\&\at{is obtained by minimizing in $c>0$}=
  \brr*{1-\tfrac{1}{2}\min\brf{\alpha,\beta}}\E{\Psi_k}
  +2\normt{\mL}\eta^2
  +2\sqrt{\tfrac{8\eta^2}{c\alpha}}
  \sqrt{\tfrac{64c\eta^2}{\alpha}\sqnt{\mL}
  +4c\alpha^2\sqnt{\mSigma}}
  \\&\at{uses the inequality $\sqrt{a+b}\leq\sqrt{a}+\sqrt{b}$}\leq
  \brr*{1-\tfrac{1}{2}\min\brf{\alpha,\beta}}\E{\Psi_k}
  +2\normt{\mL}\eta^2
  +\tfrac{32\sqrt{2}\eta^2}{\alpha}\normt{\mL}
  +8\sqrt{2\alpha}\eta\normt{\mSigma}
  \\&\at{uses the parameters in \cref{eq:parameters_cvx1}}\leq
  \brr*{1-\tfrac{1}{2}\beta}\E{\Psi_k}
  +\tfrac{1}{4}\beta\epsilon,
\end{align*}
where \annotate. Hence, we can upper-bound $\E{f(x_K) - f^*}$ as follows:
\begin{align*}
  \E{f(x_K) - f^*}
  &\at{uses the definition of $\Psi_k$ in \cref{eq:Psi1}}\leq
  \E{\Psi_K}
  \at{uses the upper bound on $\E{\Psi_k}$ above and the inequality $\sum_{k=0}^{K-1}(1-\beta/2)^k\leq2/\beta$}\leq
  \brr*{1-\tfrac{1}{2}\beta}^K\E{\Psi_0} + \tfrac{1}{2}\epsilon
  \at{uses the choice of $K$ in \cref{eq:K_cvx1}}\leq
  \epsilon,
\end{align*}
where \annotate.\qed

\newpage
\section[Proofs for \crtCref{sec:main} (Option~\ref{eq:m2})]{Proofs for \Cref{sec:main} (Option~\ref{eq:m2})}

\proofsubsection{lem:m2}

We can upper-bound $\E[\xi_{k+1}\sim \cD]{\sqn{m_{k+1} - \nabla f(x_{k+1})}_{\mB^{-1}}}$ as follows:
\begin{align*}
  \mi{1}\E[\xi_{k+1}\sim \cD]{\sqn{m_{k+1} - \nabla f(x_{k+1})}_{\mB^{-1}}}=
  \\&=
  \E[\xi_{k+1}\sim \cD]{\sqn{\E[\xi_{k+1}\sim \cD]{m_{k+1}} - \nabla f(x_{k+1})}_{\mB^{-1}}}
  +
  \E[\xi_{k+1}\sim \cD]{\sqn{m_{k+1} - \E[\xi_{k+1}\sim \cD]{m_{k+1}}}_{\mB^{-1}}}
  \\&\at{uses \cref{eq:m1} and the definition of $n(x;\xi)$ in \Cref{ass:grad}}\leq
  \sqn{(1-\alpha)m_k + \alpha \nabla f(\ox_{k+1})- \nabla f(x_{k+1})}_{\mB^{-1}}
  +
  \alpha^2\E[\xi_{k+1}\sim \cD]{\sqn{n(\ox_{k+1};\xi_{k+1})}_{\mB^{-1}}}
  \\&\at{uses \Cref{lem:B} and \Cref{ass:grad}}\leq
  \sqn{(1-\alpha)(m_k - \nabla f(x_k)) + \alpha \nabla f(\ox_{k+1}) + (1-\alpha)\nabla f(x_k) - \nabla f(x_{k+1})}_{\mB^{-1}}
  +
  4\alpha^2\sqnt{\mSigma}
  \\&\at{uses the convexity of $\sqn{}_{\mB^{-1}}$}\leq
  (1-\alpha)\sqn{m_k - \nabla f(x_k)}_{\mB^{-1}}
  +(1/\alpha)\sqn{\alpha \nabla f(\ox_{k+1}) + (1-\alpha)\nabla f(x_k) - \nabla f(x_{k+1})}_{\mB^{-1}}
  \\&
  +4\alpha^2\sqnt{\mSigma}
  \\&\at{uses the triangle inequality}\leq
  (1-\alpha)\sqn{m_k - \nabla f(x_k)}_{\mB^{-1}}
  +(3/\alpha)\sqn{\nabla^2 f(x_{k+1})(\alpha \ox_{k+1} + (1-\alpha)x_k - x_{k+1})}_{\mB^{-1}}
  \\&
  +(3(1-\alpha)^2/\alpha)\sqn{\nabla f(x_k) - \nabla f(x_{k+1}) - \nabla^2 f(x_{k+1})(x_k - x_{k+1})}_{\mB^{-1}}
  \\&
  +3\alpha\sqn{\nabla f(\ox_{k+1}) - \nabla f(x_{k+1}) - \nabla^2 f(x_{k+1})(\ox_{k+1} - x_{k+1})}_{\mB^{-1}}
  +4\alpha^2\sqnt{\mSigma}
  \\&\at{use the definition of $\ox_{k+1}$ in \cref{eq:m2}}=
  (1-\alpha)\sqn{m_k - \nabla f(x_k)}_{\mB^{-1}}
  +4\alpha^2\sqnt{\mSigma}
  \\&
  +(3(1-\alpha)^2/\alpha)\sqn{\nabla f(x_k) - \nabla f(x_{k+1}) - \nabla^2 f(x_{k+1})(x_k - x_{k+1})}_{\mB^{-1}}
  \\&
  +3\alpha\sqn{\nabla f(\ox_{k+1}) - \nabla f(x_{k+1}) - \nabla^2 f(x_{k+1})(\ox_{k+1} - x_{k+1})}_{\mB^{-1}}
  \\&\at{is implied by \Cref{ass:T} and \Cref{lem:B}}\leq
  (1-\alpha)\sqn{m_k - \nabla f(x_k)}_{\mB^{-1}}
  +\tfrac{12(1-\alpha)^2}{\alpha}\sqnt{\mT}\norm{x_k - x_{k+1}}_{\mB}^4
  +12\alpha\sqnt{\mT}\norm{\ox_{k+1} - x_{k+1}}_{\mB}^4
  \\&
  +4\alpha^2\sqnt{\mSigma}
  \\&\at{use the definition of $\ox_{k+1}$ in \cref{eq:m2}}=
  (1-\alpha)\sqn{m_k - \nabla f(x_k)}_{\mB^{-1}}
  +\tfrac{12(1-\alpha)^2}{\alpha}\brs*{1 + \tfrac{(1-\alpha)^2}{\alpha^2}}\sqnt{\mT}\norm{x_k - x_{k+1}}_{\mB}^4
  +4\alpha^2\sqnt{\mSigma}
  \\&\at{uses \Cref{lem:HR}, the definition of $\mB$ in \cref{eq:B}, and \Cref{lem:x}}\leq
  (1-\alpha)\sqn{m_k - \nabla f(x_k)}_{\mB^{-1}}
  +192\sqnt{\mT}\eta^4/\alpha^3
  +4\alpha^2\sqnt{\mSigma},
\end{align*}
where \annotate. It remains to take the full expectation.\qed

\proofsubsection{cor:m2}
We can upper-bound $\E{\sqn{m_{k+1} - \nabla f(x_{k+1})}_{\mB^{-1}}}$ as follows:
\begin{align*}
  \mi{3}\E*{\sqn{m_{k+1} - \nabla f(x_{k+1})}_{\mB^{-1}}}\leq
  \\&\at{uses \Cref{lem:m2}}\leq
  (1-\alpha)^{k+1}\E*{\sqn{m_0 - \nabla f(x_0)}_{\mB^{-1}}}
  +\brs*{192\sqnt{\mT}\eta^4/\alpha^3+4\alpha^2\sqnt{\mSigma}} \tsum_{i = 0}^k (1 - \alpha)^i
  \\&\at{uses the definition $m_0 = g(x_0; \xi_0)$ on \cref{line:init} and the definition of $n(x;\xi)$ in \Cref{ass:grad}}\leq
  (1-\alpha)^{k+1}\E*{\sqn{n(x_0;\xi_0)}_{\mB^{-1}}}
  +\brs*{192\sqnt{\mT}\eta^4/\alpha^3+4\alpha^2\sqnt{\mSigma}} \tsum_{i = 0}^k (1 - \alpha)^i
  \\&\at{uses \Cref{lem:B} and \Cref{ass:grad}}\leq
  4(1-\alpha)^{k+1}\sqnt{\mSigma}
  +\brs*{192\sqnt{\mT}\eta^4/\alpha^3+4\alpha^2\sqnt{\mSigma}} \tsum_{i = 0}^k (1 - \alpha)^i
  \\&\at{uses the inequality $\tsum_{i = 0}^k (1 - \alpha)^i \leq 1/\alpha$}\leq
  4\brs*{(1-\alpha)^{k+1} + \alpha}\sqnt{\mSigma}
  +192(\eta/\alpha)^4\sqnt{\mT}
\end{align*}
where \annotate.
Now, we can upper-bound $\E*{\norm{m_k - \nabla f(x_k)}_{\mB^{-1}}}$ as follows:
\begin{align*}
  \E*{\norm{m_k - \nabla f(x_k)}_{\mB^{-1}}}
  &\at{uses Jensen's inequality for expectations}\leq
  \sqrt{\E*{\sqn{m_k - \nabla f(x_k)}_{\mB^{-1}}}}
  \\&\at{uses the upper bound above}\leq
  \sqrt{4\brs*{(1-\alpha)^k + \alpha}\sqnt{\mSigma}
  +192(\eta/\alpha)^4\sqnt{\mT}}
  \\&\at{uses the inequality $\sqrt{a+b}\leq\sqrt{a}+\sqrt{b}$ and the inequality $\sqrt{1-\alpha} \leq 1-\alpha/2$}\leq
  2\brs*{(1-\alpha/2)^k + \sqrt{\alpha}}\normt{\mSigma}
  +8\sqrt{3}(\eta/\alpha)^2\normt{\mT},
\end{align*}
where \annotate.\qed

\proofsubsection{thm:m2}

We can upper-bound $\E*{\min_{k\in\brf{0,K-1}} \grad{\cR}{x_k}}$ as follows:
\begin{align*}
  \mi{3}
  \E*{\tmin_{k\in\brf{0,K-1}} \grad{\cR}{x_k}}\leq
  \\&\leq
  \tfrac{1}{K}\tsum_{k=0}^{K-1}\E*{ \grad{\cR}{x_k}}
  \\&\at{uses \Cref{lem:tr} and the definition of $\eta$ in \cref{eq:parameters2}}\leq
  \tfrac{1}{\eta K}\tsum_{k=0}^{K-1}\E*{ f(x_k) - f(x_{k+1}) + 4\eta\norm{\nabla f(x_k) - m_k}_{\mB^{-1}} + 2\normt{\mL}\eta^2}
  \\&\at{uses the definition of $\Delta_0$}\leq
  \tfrac{\Delta_0}{\eta K}
  +2\normt{\mL}\eta
  +\tfrac{4}{K}\tsum_{k=0}^{K-1}\E*{\norm{\nabla f(x_k) - m_k}_{\mB^{-1}}}
  \\&\at{uses \Cref{cor:m2}}\leq
  \tfrac{\Delta_0}{\eta K}
  +2\normt{\mL}\eta
  +\tfrac{4}{K}\tsum_{k=0}^{K-1}\brs*{ 2\brs*{(1-\alpha/2)^k + \sqrt{\alpha}}\normt{\mSigma}
  +8\sqrt{3}(\eta/\alpha)^2\normt{\mT}}
  \\&\at{uses the inequality $\tsum_{k=0}^{K-1}(1-\alpha/2)^k \leq 2/\alpha$}=
  \tfrac{\Delta_0}{\eta K}
  +\tfrac{16\normt{\mSigma}}{\alpha K}
  +\tfrac{32\sqrt{3}\normt{\mT}\eta^2}{\alpha^2}
  +2\normt{\mL}\eta
  +8\sqrt{\alpha}\normt{\mSigma},
\end{align*}
where \annotate. It remains to use the definitons of the parameters in \cref{eq:parameters2,eq:K2}.

Furthermore, in the convex case, we define the following Lyapunov function:
\begin{equation}\label{eq:Psi2}
  \Psi_k = f(x_k) - f^* + c\sqn{m_k - \nabla f(x_k)},
\end{equation}
where $c > 0$. We can upper-bound $\E{\Psi_{k+1}}$ as follows:
\begin{align*}
  \E{\Psi_{k+1}}
  &\at{uses \Cref{lem:tr,lem:m2}}\leq
  \E{f(x_k) - f^*} - \eta \E{\grad{\eta/\beta}{x_k}} + 4\eta\E{\norm{\nabla f(x_k) - m_k}_{\mB^{-1}}} + 2\normt{\mL}\eta^2
  \\&
  +c(1-\alpha)\E{\sqn{m_k - \nabla f(x_k)}_{\mB^{-1}}}
  +192c\sqnt{\mT}\eta^4/\alpha^3
  +4c\alpha^2\sqnt{\mSigma}
  \\&\at{uses the definition of $\eta$ in \cref{eq:parameters2} and \Cref{lem:convex}}\leq
  (1-\beta)\E{f(x_k) - f^*}
  +c(1-\alpha)\E{\sqn{m_k - \nabla f(x_k)}_{\mB^{-1}}}
  \\&
  +4\eta\E{\norm{\nabla f(x_k) - m_k}_{\mB^{-1}}}
  +2\normt{\mL}\eta^2
  +192c\sqnt{\mT}\eta^4/\alpha^3
  +4c\alpha^2\sqnt{\mSigma}
  \\&\at{uses Young's inequality}\leq
  (1-\beta)\E{f(x_k) - f^*}
  +c(1-\alpha/2)\E{\sqn{m_k - \nabla f(x_k)}_{\mB^{-1}}}
  +\tfrac{8\eta^2}{c\alpha}
  \\&
  +2\normt{\mL}\eta^2
  +192c\sqnt{\mT}\eta^4/\alpha^3
  +4c\alpha^2\sqnt{\mSigma}
  \\&\at{uses the definition of $\Psi_k$ in \cref{eq:Psi2}}\leq
  \brr*{1-\tfrac{1}{2}\min\brf{\alpha,\beta}}\E{\Psi_k}
  +2\normt{\mL}\eta^2
  +192c\sqnt{\mT}\eta^4/\alpha^3
  +4c\alpha^2\sqnt{\mSigma}
  +\tfrac{8\eta^2}{c\alpha}
  \\&\at{is obtained by minimizing in $c>0$}=
  \brr*{1-\tfrac{1}{2}\min\brf{\alpha,\beta}}\E{\Psi_k}
  +2\normt{\mL}\eta^2
  +2\sqrt{\tfrac{8\eta^2}{c\alpha}}\sqrt{192c\sqnt{\mT}\eta^4/\alpha^3 + 4c\alpha^2\sqnt{\mSigma}}
  \\&\at{uses the inequality $\sqrt{a+b}\leq\sqrt{a}+\sqrt{b}$}\leq
  \brr*{1-\tfrac{1}{2}\min\brf{\alpha,\beta}}\E{\Psi_k}
  +2\normt{\mL}\eta^2
  +\tfrac{32\sqrt{6}\eta^3}{\alpha^2}\normt{\mT}
  +8\sqrt{2\alpha}\eta\normt{\mSigma}
  \\&\at{uses the parameters in \cref{eq:parameters_cvx2}}\leq
  \brr*{1-\tfrac{1}{2}\beta}\E{\Psi_k}
  +\tfrac{1}{4}\beta\epsilon,
\end{align*}
where \annotate. Hence, we can upper-bound $\E{f(x_K) - f^*}$ as follows:
\begin{align*}
  \E{f(x_K) - f^*}
  &\at{uses the definition of $\Psi_k$ in \cref{eq:Psi2}}\leq
  \E{\Psi_K}
  \at{uses the upper bound on $\E{\Psi_k}$ above and the inequality $\sum_{k=0}^{K-1}(1-\beta/2)^k\leq2/\beta$}\leq
  \brr*{1-\tfrac{1}{2}\beta}^K\E{\Psi_0} + \tfrac{1}{2}\epsilon
  \at{uses the choice of $K$ in \cref{eq:K_cvx2}}\leq
  \epsilon,
\end{align*}
where \annotate.\qed

\newpage
\section[Proofs for \crtCref{sec:main} (Option~\ref{eq:m3})]{Proofs for \Cref{sec:main} (Option~\ref{eq:m3})}

\proofsubsection{lem:m3}

We can upper-bound $\E[\xi_{k+1}\sim \cD]{\sqn{m_{k+1} - \nabla f(x_{k+1})}_{\mB^{-1}}}$ as follows:
\begin{align*}
  \mi{0.2}\E[\xi_{k+1}\sim \cD]{\sqn{m_{k+1} - \nabla f(x_{k+1})}_{\mB^{-1}}}=
  \\&=
  \E[\xi_{k+1}\sim \cD]{\sqn{m_{k+1} - \E[\xi_{k+1}\sim \cD]{m_{k+1}}}_{\mB^{-1}}}
  +
  \E[\xi_{k+1}\sim \cD]{\sqn{\E[\xi_{k+1}\sim \cD]{m_{k+1}} - \nabla f(x_{k+1})}_{\mB^{-1}}}
  \\&\at{uses \cref{eq:m3}}=
  \E[\xi_{k+1}\sim \cD]{\sqn{(1-\alpha)(g(x_k;\xi_{k+1}) - \nabla f(x_k)) - (g(x_{k+1};\xi_{k+1}) - \nabla f(x_{k+1}))}_{\mB^{-1}}}
  \\&
  +(1-\alpha)^2\sqn{m_k - \nabla f(x_k)}_{\mB^{-1}}
  \\&\at{uses the triangle inequality and the definition of $n(x;\xi)$ in \Cref{ass:grad}}\leq
  2(1-\alpha)^2\E[\xi_{k+1}\sim \cD]{\sqn{g(x_{k+1};\xi_{k+1}) - g(x_k;\xi_{k+1})}_{\mB^{-1}}}
  +2\alpha^2\E[\xi_{k+1}\sim \cD]{\sqn{n(x_{k+1};\xi_{k+1})}_{\mB^{-1}}}
  \\&
  +(1-\alpha)^2\sqn{m_k - \nabla f(x_k)}_{\mB^{-1}}
  \\&\at{is implied by \Cref{ass:grad,ass:M} and \Cref{lem:B}}\leq
  32(1-\alpha)^2\sqnt{\mM}\sqn{x_{k+1}-x_k}_{\mB}
  +8\alpha^2\sqnt{\mSigma}
  +(1-\alpha)^2\sqn{m_k - \nabla f(x_k)}_{\mB^{-1}}
  \\&\at{uses \Cref{lem:HR} and the definition of $\mB$ in \cref{eq:B}}\leq
  (1-\alpha)^2\sqn{m_k - \nabla f(x_k)}_{\mB^{-1}}
  +128\sqnt{\mM}\eta^2
  +8\alpha^2\sqnt{\mSigma},
\end{align*}
where \annotate. It remains to take the full expectation.\qed

\proofsubsection{cor:m3}
We can upper-bound $\E{\sqn{m_{k+1} - \nabla f(x_{k+1})}_{\mB^{-1}}}$ as follows:
\begin{align*}
  \mi{3}\E*{\sqn{m_{k+1} - \nabla f(x_{k+1})}_{\mB^{-1}}}\leq
  \\&\at{uses \Cref{lem:m3}}\leq
  (1-\alpha)^{2(k+1)}\E*{\sqn{m_0 - \nabla f(x_0)}_{\mB^{-1}}}
  +\brs*{128\sqnt{\mM}\eta^2 + 8\alpha^2\sqnt{\mSigma}} \tsum_{i = 0}^k (1 - \alpha)^{2i}
  \\&\at{uses the definition $m_0 = g(x_0; \xi_0)$ on \cref{line:init} and the definition of $n(x;\xi)$ in \Cref{ass:grad}}\leq
  (1-\alpha)^{2(k+1)}\E*{\sqn{n(x_0;\xi_0)}_{\mB^{-1}}}
  +\brs*{128\sqnt{\mM}\eta^2 + 8\alpha^2\sqnt{\mSigma}} \tsum_{i = 0}^k (1 - \alpha)^{2i}
  \\&\at{uses \Cref{lem:B} and \Cref{ass:grad}}\leq
  4(1-\alpha)^{2(k+1)}\sqnt{\mSigma}
  +\brs*{128\sqnt{\mM}\eta^2 + 8\alpha^2\sqnt{\mSigma}} \tsum_{i = 0}^k (1 - \alpha)^{2i}
  \\&\at{uses the inequality $\tsum_{i = 0}^k (1 - \alpha)^{2i} \leq 1/\alpha$}\leq
  4\brs{(1-\alpha)^{2(k+1)} + 2\alpha}\sqnt{\mSigma}
  +128\sqnt{\mM}\eta^2/\alpha,
\end{align*}
where \annotate.
Now, we can upper-bound $\E*{\norm{m_k - \nabla f(x_k)}_{\mB^{-1}}}$ as follows:
\begin{align*}
  \E*{\norm{m_k - \nabla f(x_k)}_{\mB^{-1}}}
  &\at{uses Jensen's inequality for expectations}\leq
  \sqrt{\E*{\sqn{m_k - \nabla f(x_k)}_{\mB^{-1}}}}
  \\&\at{uses the upper bound above}\leq
  \sqrt{4\brs{(1-\alpha)^{2k} + 2\alpha}\sqnt{\mSigma}
  +128\sqnt{\mM}\eta^2/\alpha}
  \\&\at{uses the inequality $\sqrt{a+b}\leq\sqrt{a}+\sqrt{b}$}\leq
  2\brs*{(1-\alpha)^k + \sqrt{2\alpha}}\normt{\mSigma}
  +8\sqrt{2/\alpha}\normt{\mM}\eta,
\end{align*}
where \annotate.\qed

\proofsubsection{thm:m3}

We can upper-bound $\E*{\min_{k\in\brf{0,K-1}} \grad{\cR}{x_k}}$ as follows:
\begin{align*}
  \mi{3}
  \E*{\tmin_{k\in\brf{0,K-1}} \grad{\cR}{x_k}}\leq
  \\&\leq
  \tfrac{1}{K}\tsum_{k=0}^{K-1}\E*{ \grad{\cR}{x_k}}
  \\&\at{uses \Cref{lem:tr} and the definition of $\eta$ in \cref{eq:parameters3}}\leq
  \tfrac{1}{\eta K}\tsum_{k=0}^{K-1}\E*{ f(x_k) - f(x_{k+1}) + 4\eta\norm{\nabla f(x_k) - m_k}_{\mB^{-1}} + 2\normt{\mL}\eta^2}
  \\&\at{uses the definition of $\Delta_0$}\leq
  \tfrac{\Delta_0}{\eta K}
  +2\normt{\mL}\eta
  +\tfrac{4}{K}\tsum_{k=0}^{K-1}\E*{\norm{\nabla f(x_k) - m_k}_{\mB^{-1}}}
  \\&\at{uses \Cref{cor:m3}}\leq
  \tfrac{\Delta_0}{\eta K}
  +2\normt{\mL}\eta
  +\tfrac{4}{K}\tsum_{k=0}^{K-1}\brs*{2\brs{(1-\alpha)^k + \sqrt{2\alpha}}\normt{\mSigma}
  +8\sqrt{2/\alpha}\normt{\mM}\eta}
  \\&\at{uses the inequality $\tsum_{k=0}^{K-1}(1-\alpha)^k \leq 1/\alpha$}=
  \tfrac{\Delta_0}{\eta K}
  +\tfrac{8\normt{\mSigma}}{\alpha K}
  +\tfrac{32\sqrt{2}\normt{\mM}\eta}{\sqrt{\alpha}}
  +2\normt{\mL}\eta
  +8\sqrt{2\alpha}\normt{\mSigma},
\end{align*}
where \annotate. It remains to use the definitons of the parameters in \cref{eq:parameters3,eq:K3}.

Furthermore, in the convex case, we define the following Lyapunov function:
\begin{equation}\label{eq:Psi3}
  \Psi_k = f(x_k) - f^* + c\sqn{m_k - \nabla f(x_k)},
\end{equation}
where $c > 0$. We can upper-bound $\E{\Psi_{k+1}}$ as follows:
\begin{align*}
  \E{\Psi_{k+1}}
  &\at{uses \Cref{lem:tr,lem:m3}}\leq
  \E{f(x_k) - f^*} - \eta \E{\grad{\eta/\beta}{x_k}} + 4\eta\E{\norm{\nabla f(x_k) - m_k}_{\mB^{-1}}} + 2\normt{\mL}\eta^2
  \\&
  +c(1-\alpha)\E{\sqn{m_k - \nabla f(x_k)}_{\mB^{-1}}}
  +128c\sqnt{\mM}\eta^2
  +8c\alpha^2\sqnt{\mSigma}
  \\&\at{uses the definition of $\eta$ in \cref{eq:parameters3} and \Cref{lem:convex}}\leq
  (1-\beta)\E{f(x_k) - f^*}
  +c(1-\alpha)\E{\sqn{m_k - \nabla f(x_k)}_{\mB^{-1}}}
  \\&
  +4\eta\E{\norm{\nabla f(x_k) - m_k}_{\mB^{-1}}}
  +2\normt{\mL}\eta^2
  +128c\sqnt{\mM}\eta^2
  +8c\alpha^2\sqnt{\mSigma}
  \\&\at{uses Young's inequality}\leq
  (1-\beta)\E{f(x_k) - f^*}
  +c(1-\alpha/2)\E{\sqn{m_k - \nabla f(x_k)}_{\mB^{-1}}}
  +\tfrac{8\eta^2}{c\alpha}
  \\&
  +2\normt{\mL}\eta^2
  +128c\sqnt{\mM}\eta^2
  +8c\alpha^2\sqnt{\mSigma}
  \\&\at{uses the definition of $\Psi_k$ in \cref{eq:Psi3}}\leq
  \brr*{1-\tfrac{1}{2}\min\brf{\alpha,\beta}}\E{\Psi_k}
  +2\normt{\mL}\eta^2
  +128c\sqnt{\mM}\eta^2
  +8c\alpha^2\sqnt{\mSigma}
  +\tfrac{8\eta^2}{c\alpha}
  \\&\at{is obtained by minimizing in $c>0$}=
  \brr*{1-\tfrac{1}{2}\min\brf{\alpha,\beta}}\E{\Psi_k}
  +2\normt{\mL}\eta^2
  +2\sqrt{\tfrac{8\eta^2}{c\alpha}}
  \sqrt{128c\sqnt{\mM}\eta^2
  +8c\alpha^2\sqnt{\mSigma}}
  \\&\at{uses the inequality $\sqrt{a+b}\leq\sqrt{a}+\sqrt{b}$}\leq
  \brr*{1-\tfrac{1}{2}\min\brf{\alpha,\beta}}\E{\Psi_k}
  +2\normt{\mL}\eta^2
  +\tfrac{64\eta^2}{\sqrt{\alpha}}
  \normt{\mM}
  +16\sqrt{\alpha}\eta\normt{\mSigma}
  \\&\at{uses the parameters in \cref{eq:parameters_cvx3}}\leq
  \brr*{1-\tfrac{1}{2}\beta}\E{\Psi_k}
  +\tfrac{1}{4}\beta\epsilon,
\end{align*}
where \annotate. Hence, we can upper-bound $\E{f(x_K) - f^*}$ as follows:
\begin{align*}
  \E{f(x_K) - f^*}
  &\at{uses the definition of $\Psi_k$ in \cref{eq:Psi3}}\leq
  \E{\Psi_K}
  \at{uses the upper bound on $\E{\Psi_k}$ above and the inequality $\sum_{k=0}^{K-1}(1-\beta/2)^k\leq2/\beta$}\leq
  \brr*{1-\tfrac{1}{2}\beta}^K\E{\Psi_0} + \tfrac{1}{2}\epsilon
  \at{uses the choice of $K$ in \cref{eq:K_cvx3}}\leq
  \epsilon,
\end{align*}
where \annotate.\qed

\end{document}